\documentstyle{amsppt}

\NoRunningHeads

\mag=1200
\loadbold

\centerline{\bf Coding Theory and Uniform Distributions}
\smallskip

\centerline{By}
\smallskip

\centerline{\bf M.~M.~Skriganov, {\rm St.~Petersburg}}

\topmatter
\rightheadtext{$\text{\rm Coding Theory and Uniform Distributions}$}

\leftheadtext{$\text{\rm M.~M.~Skriganov}$}

\def\ph{\varphi}
\def\th{\theta}
\def\ga{\gamma}
\def\BM{\Bbb M}
\def\Ga{\Gamma}
\def\cd{\partial}
\def\deg{\text{\rm deg}\,}
\def\tf{\tilde f}
\def\Tr{\text{\rm Tr}\,}
\def\tsi{\tilde\sigma}
\def\s{\subset}
\def\gB{\goth B}
\def\si{\sigma}
\def\gS{\goth S}
\def\la{\lambda}
\def\BQ{\Bbb Q}
\def\BF{\Bbb F}
\def\al{\alpha}
\def\CC{\Cal C}
\def\be{\beta}
\def\om{\omega}
\def\l{\langle}
\def\r{\rangle}
\def\Mat{\text{\rm Mat}\,}
\def\BC{\Bbb C}
\def\dl{\delta}
\def\ka{\varkappa}
\def\gE{\goth E}
\def\BN{\Bbb N}
\def\Dl{\Delta}
\def\vol{\text{\rm vol}\,}
\def\CL{\Cal L}

\def\CH{\Cal H}
\def\BR{\Bbb R}
\def\Om{\Omega}

\footnote"{}"{1991 {Mathematics Subject Classification}. 11K38,
11T71, 94B60}
\footnote"{}"{{Key words and phrases}. Combinatorial structure
of uniform distributions, non-Hamming metrics on $q$-ary codes}
\footnote"{}"{The work was partially supported by Russian Fund
for Fundamental 
Research (Project No.~99-01-00106).}
\endtopmatter

{\eightpoint {\bf Abstract.} In the present paper we introduce
and study finite 
point subsets 
of a special kind, called optimum distributions, in the
$n$-dimensional unit cube. Such distributions are closely  related
with known $(\dl,s,n)$-nets of low discrepancy.
It turns out that optimum distributions have a rich
combinatorial structure. Namely, we show that optimum
distributions can be characterized completely as maximum
distance separable codes with respect to a non-Hamming metric.
Weight spectra of such codes can be evaluated precisely. 
We also consider linear codes and distributions and study their
general properties including the duality with respect to a
suitable inner product. The corresponding generalized
MacWilliams identies for weight enumerators are briefly
discussed. Broad classes of linear maximum distance
separable codes and linear optimum distributions are explicitly
constructed in the paper by the Hermite interpolations over
finite fields. 
}

\document

\vskip1cm

\centerline{\bf 1. Introduction}

\bigskip

\TagsOnRight

The present paper deals with a combinatorial structure of point
sets uniformly distributed in the $n$-dimensional unit cube. We
introduce point sets of a special kind, called optimum
distributions, which are close related with known distributions,
so-called $(\dl,s,n)$-nets, of low discrepancy. It turns out
that optimum distributions have a rich interior combinatorial
structure, namely, they can be characterized completely as
maximum distance separable codes over finite fields with respect
to a non-Hamming metric. Moreover, it is found  that weight
spectra of such codes can be evaluated precisely. 

We also study linear codes and distributions. In this case new
good codes and distributions can be obtained from given ones by
duality with respect to a suitable inner product. The
corresponding MacWilliams identities for weight enumerators in
the non-Hamming metric are also briefly discussed in the paper. 
 Notice that
broad classes of linear maximum distance separable codes and
linear optimum distributions can be explicitly given  by interpolations
over finite fields. 

We emphasize that the specifying of a 
non-Hamming metric on the space of  distributions is of
crucial importance for our study. As soon as  the metric is
introduced, the basic concepts and methods of the modern coding
theory may be applied to uniform distributions. In the present
paper this approach is illustrated by several general results. 

It is conceivable that the theory of uniform distributions may
also work for 
error-correcting codes. Further studies of intimate linkages
between coding theory and uniform distributions should be
exceptionally intriguing. 

A non Hamming metric cited above has recently been proposed by
Rosenbloom and Tsfasman [25] in coding theory. However, in the
context of  $(\dl,s,n)$-nets a
similar concept was implicitly introduced in the 1987 year
paper [20] by Niederreiter. In fact, a respective non-Hamming
weight of a linear code was defined in [20, Definitions 6.8 and
7.1] in terms of the code parity-check matrix. 

In other contexts, the metric was independently rediscovered by
several authors (see Martin and Stinson [18], and Skriganov
[28]). It seems reasonable to say that subjects of such a kind
are of great concern in many areas of combinatorial mathematics.

Since concepts and results from rather different fields are used
in the present work, necessary comments and references will be
given wherever possible. 
The reader interested in the fundamental aspects of  coding
theory is refered to the encyclopaedic book MacWilliams and
Sloane [16] and references therein.

We recall basic facts in the theory of uniform distributions. For
necessary details  we refer to Beck and Chen [4], Kuipers and
Niederreiter [12], and a recent book Matou\v sek [19]. 

Let $D\subset U^n$ be a distribution of finitely many points in
the $n$-dimensional unit cube $U^n=[0,1)^n$. The
$L_\infty$-discrepancy $\Cal L[D]$ is defined by
$$
\CL[D]=\sup_{Y\in U^n}| \#\{D\cap B(Y)\}
-\#\{D\}\vol B(Y)|,
\tag1.1
$$
where $B(Y)=[0,y_1)\ldots[0,y_n)$, $Y=(y_1,\ldots,y_n)\in U^n$,
are rectangular boxes in $U^n$, and $\vol B(Y)=y_1\ldots y_n$ . 

A main problem in the theory of uniform distributions is the
construction of distributions with minimal discrepancies (1.1).
It is known that for each $N\ge2$ there exist distributions
$D_N\subset U^n$ consisting of $N$ points with 
$$
\CL[D_N]=O((\log N)^{n-1}),
\tag1.2
$$
where the implied constant depends only on $n$. 

By a theorem of Roth for any distribution $D_N$ of $N\ge 2$
points one has the following lower bound 
$$
\CL[D_N]>c_n(\log N)^{\frac12(n-1)},
\tag1.3
$$
and by a theorem of Schmidt in two dimensions the bound (1.3)
admits the following improvements 
$$
\CL[D_N]>c\,\log N.
\tag1.4
$$
In (1.3) and (1.4) $c_n$ and $c$ are positive absolute  constants
independent of $N$. 

Thus, in two dimensions  bound (1.2) is best possible. In
dimensions $n>2$ an exact minimal order of  discrepancy (1.1)
still remains inaccessible. By [4, Sections 1.1 and 9.2] this is
a ``Great open problem'' in theory of uniform distributions. See
also Baker [2] and Beck [3] where certain improvements of the
lower bound (1.3) for  discrepancy (1.1) are given in
dimensions $n>2$. 

One of the known approachs to the construction of distributions of
low discrepancy can be outlined as follows. Fix a prime $p$ and
a prime power $q=p^e$, $e\in\BN$. We write $\BN$ for the set of
natural numbers and $\BN_0$ for the set of all non-negative
integers. Introduce elementary rectangular boxes $\Dl^M_A\subset
U^n$ by setting 
$$
\Dl^M_A=
\left[\frac{m_1}{q^{a_1}},\frac{m_1+1}{q^{a_1}}\right)\ldots
\left[
\frac{m_n}{q^{a_n}},\frac{m_n+1}{q^{a_n}}\right),
\tag1.5
$$
where $A=(a_1,\ldots,a_n)$, $M=(m_1,\ldots, m_n)\in\BN^n_0$, and
$m_j\in \{0,1,\ldots,q^{a_j}-1\}$, $1\le j\le n$. We write
$\BN^n_0$ for the 
set of all vectors in $\BR^n$ with non-negative integer
coordinates. Notice that boxes (1.5) have volume 
$$
\vol \Dl^M_A=q^{-a_1-\ldots-a_n}.
\tag1.6
$$
\medskip
{\it Definition 1.1.}  Given integers $0\le\dl\le s$. A
subset $D\subset U^n$ consisting of $q^s$ points is called a
$(\dl,s,n)$-{\it net of deficiency} $\dl$ {\it in base} $q$ if
each elementary box $\Dl^M_A$ of volume $q^{\dl-s}$ contains
exactly $q^\dl$ points of $D$. 

\medskip
Obviously $(\dl,s,n)$-nets fill out the unit cube $U^n$ very
uniformly. It is easily shown (cf. [4, Section 3.2] and [20,
Section 3]) that  discrepancy $\CL[D]$ of a $(\dl,s,n)$-net
$D$ satisfies the bound 
$$
\CL[D]=O(q^\dl s^{n-1})=O(q^\dl(\log N)^{n-1}),
\tag1.7
$$
where $N=\#\{D\}=q^s$, and the implied constants in (1.7) depend
only on $n$ and $q$. Hence, $(\dl,s,n)$-nets of bounded
deficiency $\dl$ satisfy  relation (1.2). 

We mention the following simple corollary of Definition 1.1 (cf.
[23, Lemma~9]) which shows that the base $p^e$ of a given net
can be always reduced to $p$ at the cost of an increase in
deficiency. 

\medskip
{\bf Lemma 1.1.} {\it Every $(\dl,s,n)$-net in base $p^e$ is a
$(\dl',es,n)$-net of deficiency $\dl'=e\dl+(e-1)(n-1)$ in base
$p$. 
}
\medskip

The first constructions of $(\dl,s,n)$-nets were given by Sobol
[29, 30] and Faure [11]. Namely, in an arbitrary dimension $n$
for any $s\in\BN$ $(\dl,s,n)$-nets of deficiency $\dl=O(n\log
n)$ in base 2 were constructed in [29, 30], and $(0,s,n)$-nets
of  zeroth deficiency in arbitrary prime base $p\ge n-1$ were given
in [11]. 

A systematic study of $(\dl, s,n)$-nets (including nets in
arbitrary integer bases $b\in\BN)$ was started by Niederreiter
[20], see also a book [21] and a recent survey [22]. Among other
things the following impressive result was recently obtained
within the context of these studies. Niederreiter and Xing [23,
24] discovered that very good $(\dl,s,n)$-nets of deficiency
$\dl=O(n)$ can be constructed in terms of algebraic curves over
finite fields. Moreover, for a given base $q$ the bound
$\dl=O(n)$ turns out to be best possible as $n\to\infty$. 

It was observed in Adams and Shader [1], Clayman and Mullen [8],
Edel and Bierbrauer [10], and Lawrence {\it et al.} [13] that
properties of $(\dl,s,n)$-nets can be improved in some instances
by error-correcting codes with large Hamining weights. The
corresponding details and further references can also be found
in the survey [24, Section 6]. 

It should be noted that $(\dl,s,n)$-nets can be used to
construct new remarkable classes of uniform distributions. For
example, probabilistic variations of $(\dl,s,n)$-nets were used
in Chen [6] (see also [4, Section 3.4]) to determine exact
minimal orders of the $L_W$-discrepancies for all $1< W<\infty$.

Concerning alternative constructions of distributions of low
discrepancy, we should mention the author's paper [26] (see also
[27]) where methods of the geometry of numbers had been used. In
a sense the approach given in the present paper is based on
similar geometric ideas adapted for finite fields. 

Main concepts of the geometry of numbers (like sphere packings,
homogeneous and inhomogeneous minima, etc.) arranged for finite
fields can be found in Conway and Sloane [9]. For general facts
on geometry of vector spaces over finite fields we refer to Lang
[14] and Lidl and Niederreiter [15]. 

In the present paper we study point sets of a special kind which
are defined as follows. Introduce the following collection of
elementary boxes (1.5) 
$$
\gE_s(q,n)=\{
\Dl^M_A:A=(a_1,\ldots,a_n),\,0\le a_j\le s,\,1\le j\le n\}.
\tag1.8
$$

\medskip
{\it Definition 1.2.} {Given an integer $0\le k\le n$. A
subset $D\subset U^n$ consisting of $q^k$ points is called an
{\it optimum} $[ns,k]_s$-{\it distribution in base} $q$ if
each elementary box $\Dl^ M_A\in\gE_s(q,n)$ of volume $q^{-k}$
contains exactly one point of $D$. 
}
\medskip

{\it Remark 1.1.}
A meaning of Definition 1.2 is to settle on a relation between
the number of points in a given distribution and the magnitude
of denominators in the definition of elementary boxes (1.5). 

It should be pointed out that in Definition~1.2 we may assume,
without loss of generality, that $s\le k\le ns$. Indeed, for
$0\le k<s$ each elementary box $\Dl^M_A\in\goth E_s(q,n)$ of
volume $q^{-k}$ also belongs to a subcollection $\goth
E_k(q,n)$, and the corresponding point set turns out to be an
optimum $[nk,k]_k$-distribution. At the same time, for $s<k\le
ns$ such rescaling is impossible. 

In the sequel an arbitrary subset $D\subset U^n$ consisting of
$q^k$ points is conveniently  said to be an $[ns,k]_s$-{\it
distribution in base} $q$ (with respect to a collection of
elementary boxes (1.8)).

The following result is a simple corollary of Definition~1.2.

\medskip
{\bf Lemma 1.2.} {\it Let $D\subset U^n$ be an optimum $[ns,k]_s$
-- distribution in base $q$. Then,

{\rm (i)} Each elementary box $\Dl^M_A\in\gE_s(q,n)$ contains
exactly $q^{k-a_1-\ldots-a_n}$ points of $D$, provided that
$a_1+\ldots+a_n\le k$. 

{\rm (ii)} Each elementary box $\Dl^M_A\in\gE_s(q,n)$ contains
at most one point of $D$, provided that $a_1+\ldots +a_n>k$. 
}
\medskip

{\it Proof.} (i) Choose integers $c_1,\ldots,c_n$ to satisfy the
relations $c_1+\ldots+c_n=k-a_1-\ldots-a_n$, $0\le c_j\le s$,
$1\le j\le n$. Let $b_j=a_j+c_j$, $1\le j\le n$, then from (1.5)
we obtain that 
$$
\Dl^M_A=\left[
\frac{q^{c_1}m_1}{q^{b_1}},\frac{q^{c_1}m_1+q^{c_1}}{q^{b_1}}\right)
\ldots 
\left[
\frac{q^{c_n}m_n}{q^{b_n}},\frac{q^{c_n}m_n+q^{c_n}}{q^{b_n}}\right).
$$
Hence, the box $\Dl^M_A$ is a union of $q^{c_1+\ldots+c_n}$
 disjoint elementary boxes
$\Dl^T_B\in\gE_s(q,n)$ of volume $q^{-k}$, and (i) follows
from Definition~1.2.

(ii) It suffices to notice that each elementary box
$\Dl^M_A\in\gE_s(q,n)$ with $a_1+\ldots+a_n>k$ is contained in
a bigger elementary box $\Dl^T_B\in\gE_s(q,n)$ of volume
$q^{-k}$, and (ii) follows from Definition~1.2. $\square$
\medskip

Now a linkage between optimum distributions and nets can be
described as follows. 

\medskip
{\bf Proposition 1.1.} {\it An optimum $[ns,k]_s$-distribution
with $s\le k\le ns$ is a $(k-s,k,n)$-net in the same base $q$.
}
\medskip

{\it Proof.} It suffices to compare Lemma 1.2 (i) with
Definition~1.1 and take into account that any elementary box
$\Dl^M_A$ of volume $q^{-s}$ belongs to the collection
$\gE_s(q,n)$. $\square$
\medskip

From Proposition~1.1 and bound (1.7) we conclude that the
discrepancy $\CL[D]$ of an optimum $[ns,k]_s$-distribution $D$
with $s\le k\le ns$ 
satisfies the bound 
$$
\CL[D]=O(q^{k-s}s^{n-1})=
O(N^{1-\frac sk}(\log N)^{n-1}),
\tag1.9
$$
where $N=\#\{D\}=q^k$. Thus, optimum $[ns,s]_s$ --
distributions coinside with $(0,s,n)$-nets of zeroth deficiency
and they satisfy  relation (1.2). 

At first glance the foregoing implies that only optimum
distributions with $k=s$ are of interest to the theory of uniform
distributions. However, such is not the case, and distributions
with $k>s$ can also be used to construct new ones which already
satisfy  relation (1.2). More precisely, from a given optimum
$[gns,gs]_s$-distribution $D\subset U^{gn}$, $g\in\BN$, we can
obtain a new optimum $[gns,gs]_{gs}$-distribution
$\pi_{g,n}D\subset U^n$ by a mapping $\pi_{g,n}:U^{gn}\to U^n$
which is merely the known Peano's bijection between points of
the unit cubes $U^{gn}$ and $U^n$.  Moreover, it turns out that
the distributions 
$\pi_{g,n}D$ have additional intriguing properties.  Concerning
Peano's mappings, the reader is refered to any text-book on 
elementary set theory. 

The study of optimum $[ns,k]_s$-distributions is a main goal
of the present paper. Previously a part of our results was given
in the author's paper [28]. 

The paper is organized as follows. Section 2 contains general
results on distributions and their $q$-ary codes. Our main
concepts of a non-Hamming metric $\rho$ and maximum distance
separable (or briefly MDS) codes in the metric $\rho$ are
introduced in Sec.~2.2 and Sec.~2.3, respectively. 

In Sec.~2.3 optimum distributions are characterized completely
in terms of their $q$-ary codes. We show that a given
distribution is an optimum distribution if and only if its code
is an MDS code in the metric $\rho$ (see Theorem~2.1).

Notice that by definition MDS codes consists of the most widely
spaced code words. For the Hamming metric $\varkappa$ a general
theory of MDS codes is given in [16, Chapter~11]. Metrics
different from the Hamming one, say, the Lee metric, are known
in coding theory and its applications (cf. [16]). Not all of
them enable a rich theory of MDS codes. In Sec.~2.5 we shall
mention a class of metrics, containing metrics $\varkappa$ and
$\rho$, which could be of interest from such standpoint (see
Remark ~2.1).  Notice also that a nontrivial group of linear
transformations preserving the metric $\rho$ can be given
explicitly (see Remark~2.2). 

It is an interesting question whether the non-Hamming metric
$\rho$ can be practically applied to concrete communication
systems. Some examples can be found in [25]. We shall also give
examples of such a kind (see Remark~2.3) in order to illustrate
main distinctions in applications of the metric $\rho$ to coding
theory and uniform distributions. 

Section 3 contains explicit formulas for weight spectra of
optimum distributions and MDS codes in the metric $\rho$ (see
Theorem~3.1). We emphasize that such formulas are fully inspired
by coding theory where similar results are well known for
MDS codes in the Hamming metric (cf. [16, Chapter 11, Theorem
6]). 

Explicit formulas for weight spectra are proved in
Section~3 for arbitrary (linear or non-linear) MDS codes and
optimum distributions. Furthermore, such formulas imply certain
necessary conditions for the existence of MDS codes 
and optimum distributions (see Proposition~3.1).

Section 4 contains results on linear codes and distributions. An
additional group structure of these subjects allows the
application of the
Fourier analysis on Abelian groups to study linear codes and
distributions. Here this approach is illustrated by several
results of a general character. 

We define the duality with respect to an inner product on the
space of codes and distributions (see Definition~4.2), and show
that the dual subjects to linear MDS $[ns,k]_s$-codes and optimum
$[ns,k]_s$-distributions are linear MDS $[ns,ns-k]_s$-codes, and
correspondingly, optimum $[ns,ns-k]_s$-distributions (see
Theorem~4.1). 

In particular, linear optimum $[ns,s]_s$- and
$[ns,(n-1)s]_s$-distributions are mutually dual. By the way, this
confirms once more that the consideration of optimum
$[ns,k]_s$-distributions with $k>s$ is of interest to us. 

It turns out that linear $(\dl,s,n)$-nets can also be
characterized completely in terms of weights of their dual
subjects (see Theorem~4.2). Such results on $(\dl,s,n)$-nets can
be treated as a metric interpretation of those given earlier in
[20]. 

Furthermore, in Sec.~4.4 we shall take a quick look at
generalized MacWilliams relations for weight enumerators in the
metric $\rho$. However, an extended discussion of this matter is
beyond the scope of the paper. 

Notice that an alternative approach to linear $(\dl,s,n)$-nets,
including MacWil\-li\-ams type theorems, has  recently been developed by
Martin and Stinson [18] in terms of association schemes. 

In Sec.~4.5 we evaluate the non-Hamming weight of a linear code
in terms of its parity-check matrix (see Proposition~4.2). This
simple result confirms in passing that the corresponding
definitions given in [20] and [25] are equivalent for linear
codes. Certainly, the definition based on the metric $\rho$ is
preferable from a geometric standpoint.

Section 5 contains explicit examples of linear MDS codes in the
metric $\rho$ and linear optimum distributions constructed by
the Hermite interpolations over finite fields. A general theory
of the Hermite interpolation problem over the field of complex
numbers $\BC$ can be found in Berezin and Zhidkov [5, Chapter 2,
Section 11]. It should be pointed out that in the problem over
finite fields we replace the usual formal derivatives by
hyperderivatives introduced by Hasse and Teichm\"uller (see [15,
Section~6.4]
and references therein). Such an  approach turns out to be very
convenient for calculations in fields of finite characteristic. 

Constructions of MDS codes in the Hamming metric by
interpolations are also known in coding theory (see, for
example, 
Mandelbaum [17]). From this standpoint our constructions of MDS
codes in the metric $\rho$ can be regarded as a suitable
generalization of the known Reed--Solomon codes, or more
generally, the redundant residue codes (see [16, Chapter 10]).
Notice that interpolations with usual formal derivatives were
also used in [25].

It was mentioned above that new algebraic-geometric
constructions of very good $(\dl,s,n)$-nets had been given by
Niderreiter and Xing [23, 24]. In the spirit of our discussion
of parallels between coding theory and uniform distributions
these constructions should be interpreted as the known Goppa
codes adapted for the non-Hamming metric $\rho$. Concerning the
Goppa codes and related topics we refer to [16, Chapter~12],
Stichtenoth [31], and Tsfasman and Vl\v adut [32]. Notice also
that algebraic-geometric codes for the metric $\rho$ are briefly
discussed in [25]. 
Unfortunately, in the present paper we have to leave aside these
extremely interesting questions. 

Section 6 contains more specific results. Here we study the
behavior of distributions under Peano's bijection
$\pi_{g,n}:U^{gn}\to U^n$. It is worth noting that the
reconstruction of codes and distributions by the mapping
$\pi_{g,n}$ is also inspired by coding theory, where
reconstructions of a similar kind are known, say, for the
so-called concatenated codes (cf. [16, Chapter 10]). 

We show that for a given optimum $[gns,gk]_s$-distribution
$D\subset U^{gn}$, $g\in\BN$, the image $\pi_{g,n}D\subset U^n$
is an optimum $[ngs,gk]_{gs}$-distributions. 
It turns out that the distribution $\pi_{g,n}D$ also has a large
weight in the usual Hamming metric $\ka$ (see Propositions~6.1
and 6.3). We give explicit examples of codes and distributions
with large weights simultaneously in both metrics $\rho$ and
$\ka$ (see Theorems~6.1 and 6.2). These topics are close
related with a recent joint work of Chen and the author [7].
Conceivably, the treatment of codes and distributions
simultaneously in several different (or even in all possible)
metrics should be of much interest to clarify the structure of
these subjects.

The author hopes that the foregoing discussion  could show the
benefits of 
interactions between coding theory and the theory of uniform
distributions. 
\bigskip

{\it Acknowledgments.} I would like to express my sincere thanks
to W.~W.~L.~Chen, S.~T.~Do\-u\-g\-he\-r\-ty, J.~Matou\v sek, and
H.~Niederreiter for numerous discussions and valuable
suggestions on the matter of the present paper.
 
\medskip
\centerline{\bf 2. Distributions and Codes}
\medskip

{\bf 2.1 Preliminaries.}
Consider the $q$-ary representation of a point $x\in[0,1)$, 
$$
x=\sum_{i\ge1} \eta_i(x)q^{-i},\quad
\eta_i(x)\in\{0,1,\ldots,q-1\}.
\tag2.1
$$
Representation (2.1) is uniquely defined if we agree that
the series (2.1) is finite for rational $q$-ary points $x=\frac
m{q^a}$, $a\in\BN$, $m\in\{0,1,\ldots,q^a-1\}$. We write
$\BQ(q^s)$, $s\in\BN$, for the set of all points $x=\frac
m{q^a}\in[0,1)$ with  $a\le s$, and $\BQ^n(q^s)$ for the set of
all points $X=(x_1,\ldots,x_n)^T\in U^n$ with coordinates
$x_j\in\BQ(q^s)$, $1\le j\le s$. Here $(\boldkey.)^T$ denotes the
transpose of a matrix $(\boldkey.)$.

For each $x\in[0,1)$ we define a projection $\tau_sx$ onto
$\BQ(q^s)$ by
$$
\tau_s x=\sum^s_{i=1}\eta_i(x)q^{-i}.
\tag2.2
$$
Obviously, we have 
$$
x-\tau_s x=\sum_{i\ge s+1}\eta_i(x)q^{-i}\in [0,q^{-s}).
\tag2.3
$$

Similarly, for each point $X=(x_1,\ldots,x_n)^T\in U^n$ we
define a projection $\tau_sX$ onto $\BQ^n(q^s)$ by 
$$
\tau_sX=(\tau_s x_1,\ldots,\tau_sx_n)^T.
\tag2.4
$$

\medskip{\bf Lemma 2.1.} {\it Let $D\subset U^n$ be an optimum
$[ns,k]_s$-distribution, then its projection $\tau_sD\subset
\BQ^n(q^s)$ is also an optimum $[ns,k]_s$-distribution.}
\medskip

{\it Proof.} From  relations (2.2), (2.3), and (2.4) we
conclude that a point $X$ falls into an elementary box
$\Dl^M_A\in \gE_s(q,n)$ if and only if its projection $\tau_sX$
does, and Lemma~2.1 follows at once from Definition~2.1.
$\square$ 
\medskip

Thus, without loss of generality, we may consider in the sequel
only optimum $[ns,k]_s$-dis\-tri\-bu\-ti\-on which are subsets in
$\BQ^n(q^s)$. 

Representation (2.1) is convenient to write in the form
$$
x=\sum^s_{i=1}\xi_i(x)q^{i-s-1},
\tag2.5
$$
where $\xi_i(x)=\eta_{s+1-i}(x)$. 

Coefficients $\xi_i(x)$ in (2.5) can be identified with
elements of a finite field. Let $\BF_q$ be a finite field with
$q=p^e$ elements. Fix a basis in $\BF_q$ over $\BF_p$, and
represent elements $\mu=(\mu^{(1)},\ldots,\mu^{(e)})\in\BF_q$ as
$e$-term sets of elements $\mu^{(i)}\in\BF_p$. In so doing,
elements of the field $\BF_p$ are thought of as reduced residues
$\{0,1,\ldots,p-1\}$ modulo $p$ (see [15] for details). With
these conventions we have the following bijection between
elements $\mu\in\BF_q$ and integers $m\in\{0,1,\ldots,q-1\}$: 
$$
\BF_q\ni \mu=(\mu^{(1)},\ldots,\mu^{(e)})\rightleftarrows m
=\sum^e_{i=1}\mu^{(i)}p^{i-1}\in\{0,1,\ldots,q-1\}.
\tag2.6
$$

By the use of bijection (2.6), coefficients $\xi_i(x)$
in (2.5) may be regarded as elements of the field $\BF_q$, and
therefore, a point $x\in\BQ(q^s)$ can be identified with the
following row matrix 
$$
\om\l x\r =\om_q \l x\r =(\xi_1(x),\ldots,\xi_s(x))\in
\Mat_{1,s}(\BF_q). 
\tag2.7
$$
We write $\Mat_{n,s}(k)$ for the linear space of all matrices
with $n$ rows and $s$ columns over a field $k$. 

Using (2.7), we identify each point $X\in\BQ^n(q^s)$ with the
following matrix, the {\it code word} of $X$,
$$
\Om\l X\r =(\om(x_1),\ldots,\om(x_n))^T\in\Mat_{n,s}(\BF_q).
\tag2.8
$$

For a given distribution $D\subset \BQ^n(q^s)$ we define its
code $\CC\l D\r\subset \Mat_{n,s}(\BF_q)$ as the following array
of the corresponding code words (2.8): 
$$
\CC\l D\r=\{\Om(X),X\in D\}.
\tag2.9
$$
Obviously, each subset $\CC\subset\Mat_{n,s}(\BF_q)$ is  code
(2.8) for a distribution $D=D\l \CC\r\subset \BQ^n(q^s)$. 

By  relations (2.5), (2.7), and (2.8) the structure of the
vector space\break $\Mat_{n,s}(\BF_q)$ can be transfered onto the set
of all rational $q$-ary points $\BQ^n(q^s)$. Namely, for any two
points $X=(x_1,\ldots,x_n)^T$ and $Y=(y_1,\ldots,y_n)^T\in
\BQ^n(q^s)$ and two elements $\al$ and $\be\in\BF_q$ we define a
linear combination $\al X\oplus\be Y=(\al x_1\oplus \be
y_1,\ldots, \al x_n \oplus \be y_n)^T\in\BQ^n(q^s)$ by setting 
$$
\xi_i(\al x_j\oplus\be y_j)=\al\xi_i(x_j)+\be\xi_i(y_i),\quad
1\le i\le s,\quad 1\le j\le n,
\tag2.10
$$
where $\xi_i(x)\in\BF_q$ are coefficients in  representation
(2.5). 

With respect to  arithmetic operations (2.10) the set
$\BQ^n(q^s)$ is a vector space of dimension $ns$ over $\BF_q$. 

We need to describe special affine subspaces 
 in $\BQ^n(q^s)$ associated with elementary boxes
(1.5). 

\medskip{\bf Lemma 2.2.} {\it For each elementary box $\Dl^M_A\in
\gE_s(q,n)$ the intersection $V^M_A=\BQ^n(q^s)\cap\Dl^M_A$ is an
affine subspace in $\BQ^n(q^s)$ of dimension $ns-a_1-\ldots-a_n$
parallel to the subspace $V^0_A=\BQ^n(q^s)\cap\Dl^0_A$, that is,
$V^M_A=V^0_A\oplus Y^M_A$, where a vector $Y^M_A\in\BQ^n(q^s)$
is defined uniquely  up to translations along $V^0_A$; in
particular one 
may take
$Y^M_A=\left(\frac{m_1}{q^{a_1}},\ldots,\frac{m_n}{q^{a_n}}\right)^T$.

More precisely, the affine subspace $V^M_A$ consists of points
$X=(x_1,\ldots,x_n)^T \in\BQ^n(q^s)$ which satisfy the following
system of $a_1+\ldots+a_n$ linear independent equations 
$$
\xi_i(x_j)=\xi_i\left(\frac{m_j}{q^{a_j}}\right),\quad
s\ge i\ge s+1-a_j,\quad 1\le j\le n,
\tag2.11
$$
where $\xi_i(\boldkey.)$ are   coefficients in  $q$-ary
representation $(2.5)$.
}
\medskip

{\it Proof.} From representation (2.1) we easily conclude
that a point $x\in\BQ(q^s)$ falls into an interval $\left[\frac
m{q^a},\frac {m+1}{q^a}\right)$, $a\in\BN_0$, $0\le a\le s$,
$m\in\{0,1,\ldots, q^a-1\}$, if and only if
$\eta_i(x)=\eta_i\left( \frac m{q^a}\right)$ for $1\le i\le a$.
Writing these equations in terms of coefficients
$\xi_i(\boldkey.)$ in (2.5), and using  definition (1.8), we
obtain (2.11). It is evident that  equations (2.11) imply all
 other statements of Lemma~2.2. $\square$
\medskip

{\bf 2.2. Definition of a non-Hamming Metric.}
We introduce the following non-Hamming {\it weight} $\rho$ on
the space $\Mat_{n,s}(\BF_q)$ (cf. [25]).
 At first, let $n=1$ and
$\om=(\xi_1,\ldots,\xi_s)\in \Mat_{1,s}(\BF_q)$. Then we put
$\rho(0)=0$, and 
$$
\rho(\om)=\rho_q(\om)=\max\{j:\xi_j\ne0\}
\tag2.12
$$
for $\om\ne0$. 

Notice that in the sequel the base $q$ is usually assumed to be
fixed, and we drop it from our notation. Only in Sec.~6.4 we
shall need to indicate the dependence on $q$ in order to treat
the variation of $q$ from $p^e$ to $p$. 

It is obvious that $\rho(\om)>0$ for $\om\ne0$,
$\rho(\al\om)=\rho(\om)$ for $\al\in\BF^*_q=\BF_q\setminus\{0\}$, and 
$$
\rho(\om+\om')\le \max\{\rho(\om),\rho(\om')\}\le
\rho(\om)+\rho(\om')
$$
for all $\om,\om'\in\Mat_{1,s}(\BF_q)$. Moreover, we have
trivially that 
$$
\rho(\om)\le s
\tag2.13
$$
for all $\om\in\Mat_{1,s}(\BF_q)$. 

Let $\Om=(\om^{(1)},\ldots,\om^{(n)})^T\in\Mat_{n,s}(\BF_q)$, 
$\om^{(j)}\in\Mat_{1,s}(\BF_q)$, $1\le j\le n$. Then we put
$$
\rho(\Om)=\rho_q(\Om)=\sum^n_{j=1}\rho(\om^{(j)}).
\tag2.14
$$

For example, let $n=2$, $s=3$, and $\Om=\left(\matrix 1&1&0\\
0&0&1\endmatrix\right)$, then
$\rho(\Om)=\rho((1\,1\,0))+\rho((0\,0\,1))= 2+3=5$.

Obviously, $\rho(\Om)>0$ for $\Om\ne0$, $\rho(\al\Om)=\rho(\Om)$
for $\al\in\BF^*_q$, and
$$
\rho(\Om+\Om')\le\rho(\Om)+\rho(\Om')
$$
for all $\Om$, $\Om'\in\Mat_{n,s}(\BF_q)$. Thus,
$\rho(\Om-\Om')$ is a {\it metric} (or a {\it distance}) on the
space $\Mat_{n,s}(\BF_q)$. 

Recall that by definition the {\it Hamming weight}
$\ka(\Om)=\ka_q(\Om)$, $\Om\in\Mat_{n,s}(\BF_q)$ is equal to the
number of non-zeroth entries in the matrix $\Om$ (cf. [16]). It
is easily seen that the weights $\ka$ and $\rho$ are related by 
$$
\ka(\Om)\le\rho(\Om)\le s\ka(\Om),
\tag2.15
$$
and these inequalities can not be improved on the whole space
$\Mat_{n,s}(\BF_q)$. In this sense the metric $\rho$ is stronger
than $\ka$ for large $s$. 
\medskip
{\bf 2.3. MDS Codes in the Metric $\rho$.}
Given an arbitrary metric $\la$ on the space
$\Mat_{n,s}(\BF_q)$. Using  relations (2.5), (2.7), and
(2.8), the metric $\la$ can be transfered onto the set of all
rational $q$-ary points $\BQ^n(q^s)$ by setting 
$$
\la(X)=\la(\Om\l X\r)
\tag2.16
$$
for the $\la$-weight of a point $X\in\BQ^n(q^s)$, and 
$$
\la(X\ominus X')=\la(\Om\l X\r-\Om\l X'\r)
\tag2.17
$$
for the $\la$-distance between two points $X,X'\in\BQ^n(q^s)$.
In (2.17) the symbol $\ominus$ denotes the subtraction with
respect to the addition $\oplus$, which was defined in (2.10). 

For each code $\CC\subset \Mat_{n,s}(\BF_q)$, which is assumed
to contain at least two elements, its {\it weight} $\la(\CC)$ is
defined by 
$$
\la(\CC)=\min\{\la(\Om-\Om'):\Om,\Om'\in\CC,\Om\ne\Om'\} 
\tag2.18
$$

Similarly, for each distribution $D\subset\BQ^n(q^s)$, which
contains at least two points, its {\it weight} $\la(D)$ is
defined by 
$$
\la(D)=\min\{\la(X\ominus  X'):X,X'\in D,X\ne X'\}. 
\tag2.19
$$

It is obvious that $\la(D)=\la(\CC\l D\r)$ and $\la(\CC)=\la(D\l
\CC\r)$. 

In the sequel $\la$ is equal to either $\ka$ or $\rho$. So that,
we write $\ka(\CC)$ and $\rho(\CC)$, and correspondingly
$\ka(D)$ and $\rho(D)$, for the Hamming and non-Hamming weights
of codes $\CC\subset\Mat_{n,s}(\BF_q)$ and distributions
$D\subset\BQ^n(q^s)$. 

Relations (2.15) imply the similar inequality for weights 
$$
\ka(C)\le \rho(C)\le s\ka(C),\quad
\ka(D)\le\rho(D)\le s\ka(D).
\tag2.20
$$

Notice also that if a code or a distribution coincides with the
whole space, then 
$$
\rho(\Mat_{n,s}(\BF_q))=\rho(\BQ^n(q^s))=1.
\tag2.21
$$
These relations follow at once from  definitions (2.12),
(2.14), and (2.18), (2.19). 

Given an integer $0\le k\le ns$. An arbitrary code
$\CC\subset\Mat_{n,s}(\BF_q)$ consisting of $q^k$ elements is
conveniently said to be an $[ns,k]_s$-{\it code}. 

The following basic result was given in [25, Theorem 1]. For the
reader convenience we also reproduce its proof.   

\medskip{\bf Proposition 2.1.} {\it For every code $\Cal
C\s\Mat_{n,s}(\Bbb F_q)$ its cardinality $\#\{C\}$ and weight
$\rho(\Cal C)$  are related by 
$$
\#\{\CC\}\le q^{ns-\rho(\CC)+1}.
\tag2.22
$$

In particular, every $[ns,k]_s$-code satisfies the bound 
$$
\rho(\CC)\le ns-k+1.
\tag2.23
$$
}
\medskip

{\it Proof.} In a matrix
$\Om=(\xi^{(j)}_i)\in\Mat_{n,s}(\BF_q)$, $1\le j\le n$, $1\le
i\le s$, we enumerate entries $\xi^{(j)}_i$ lexicographicaly,
and mark first $\rho(\CC)-1$ positions. Now, from definition
(2.18) of the weight $\rho(\CC)$ we conclude that for two
distinct code words $\Om_1,\Om_2\in\CC$ their entries can not
coinside in all other $ns-\rho(\CC)+1$ positions, because
otherwise $\rho(\Om_1-\Om_2)\le\rho(\CC)-1$. Therefore, relation
(2.22) holds and bound (2.23) follows. $\square$

Proposition 2.1 leads immediately to the following. 

\medskip
{\it Definition 2.1.} {An $[ns,k]_s$-code $\CC\subset
\Mat_{n,s}(\BF_q)$ is called a {\it maximum distance separable}
(or briefly {MDS}) $[ns,k]_s$-{\it code in the  metric} $\rho$
if its weight $\rho(\CC)= ns-k+1$. 
}
\medskip

\medskip
{\bf 2.4. Optimum Distributions and MDS Codes in the Metric
$\rho$.} 
Our first result on optimum distributions can be stated as
follows. 

\medskip{\bf Theorem 2.1.} {\it Let $D\subset\BQ^n(q^s)$ be a
distribution and $\CC\l D\r\subset \Mat_{n,s}(\BF_q)$ be its
code $(2.9)$. Then the following two statements are equivalent: 

{\rm(i)} $D$ is an optimum $[ns,k]_s$-distribution, 

{\rm(ii)} $\CC\l D\r$ is an MDS $[ns,k]_s$-code in the metric
$\rho$. 
}
\medskip

For the proof of Theorem~2.1 we need the following. 

\medskip{\bf Lemma 2.3.} {\rm(i)} {\it A given point
$X\in\BQ^n(q^s)$ falls 
into an elementary box $\Dl^0_A\in\gE_s(q,n)$ of volume
$q^{-k}$ if and only if the weight $\rho(X)\le ns-k$. 

{\rm(ii)} Two given points $X,X'\in\BQ^n(q^s)$ fall
simultaneously into an elementary box $\Dl^M_A\in\gE_s(q,n)$ of
volume $q^{-k}$ if and only if the distance $\rho(X\ominus
X')\le ns-k$. 
}
\medskip

{\it Proof.} First of all, we notice that definition (2.11)
and  representation (2.5) imply the inequality 
$$
x=\sum^{\rho(x)}_{i=1}\xi_i(x)q^{i-s-1}\le(q-1)\sum^{\rho(x)}_{i=1}
q^{i-s-1}<q^{\rho(x)-s}
\tag2.24
$$
for all $x\in\BQ(q^s)$ and the inequiality
$$
x=\sum^{\rho(x)}_{i=1}\xi_i(x)q^{i-s-1}\ge q^{\rho(x)-s-1}
\tag2.25
$$
for $x>0$, $x\in\BQ(q^s)$. 

(i) Let $\rho(X)\le ns-k$. Writing inequality (2.24) for
all coordinates of the point $X=(x_1,\ldots,x_n)^T\in\BQ^n(q^s)$
and using  definition (1.5), we conclude that
$X\in\Dl^0_{A(X)}$ with $A(X)=(a_1(x_1),\ldots,a_n(x_n))$ where
$a_j(x_j)=s-\rho(x_j)$, $1\le j\le n$. Therefore, we have (cf.
(2.14), (2.17))
$$
a_1(x_1)+\ldots+a_n(x_n)=ns-\rho(X)\ge k,
$$
and $0\le a_j(x_j)\le s$, $1\le j\le n$. Hence (cf. (1.6),
(1.8)), there exists an elementary box $\Dl^0_A\in\gE_s(q,n)$ of
volume $q^{-k}$ such that $\Dl^0_A\supseteq A^0_{A(X)}\ni X$.  

Let $\rho(X)> ns-k$. Writing  inequality (2.25) for all
non-zeroth coordinates of the point
$X=(x_1,\ldots,x_n)^T\in\BQ^n(q^s)$ and using  definitions
(1.5) and (1.8), we observe the following: if the point $X$
falls into an elementary box $\Dl^0_A\in\gE_s(q,n)$ with
$A=(a_1,\ldots,a_n)$, then $0\le a_j\le s-\rho(x_j)$ for $x_j>0$
and $0\le a_j\le s$ for $x_j=0$, $1\le j\le n$. Since
$\rho(0)=0$, we have 
$$
a_1+\ldots+a_n\le ns-\rho(x_1)-\ldots-\rho(x_n)=ns-\rho(X)<k.
$$
Hence (cf. (1.6)), the elementary box $\Dl^0_A\in\gE_s(q,n)$ has
volume strictly bigger than $q^{-k}$. 
The proof of the statement (i) is complete.

(ii) Since the distance $\rho(X\ominus X')$ is invariant under
translations $X\to X\oplus Y$, $X'\to X'\oplus Y$,
$Y\in\BQ^n(q^s)$, the statement (ii) follows at once from the
statement (i) with the help of Lemma~2.2. $\square$
\medskip

{\it Proof of Theorem~2.1.} (i) Let $D$ be an optimum
$[ns,k]_s$-distribution. Then, by Definition 1.2, each elementary
box $\Dl^M_A\in\gE_s(q,n)$ of volume $q^{-k}$ contains exactly
one point of $D$, and by Lemma~2.3 (ii), the distance
$\rho(X\ominus X')> ns-k$ for any pair of distinct points $X$,
$X'\in D$. Therefore, by Definition~2.1, $\CC\l D\r$ is an MDS
$[ns,k]_s$-code in the metric $\rho$.

(ii) Let $\CC\l D\r$ be an MDS $[ns,k]_s$-code in the metric
$\rho$. Then, by Definition~2.1, the distance $\rho(X\ominus
X')> ns-k$ for any pair of distinct points $X,X'\in D$, and by
Lemma~2.3 (ii), each elementary box $\Dl^M_A\in\gE_s(q,n)$ of
volume $q^{-k}$ contains at most one point of $D$. Now we
observe that for a given $A=(a_1,\ldots,a_n)\in\BN^n_0$ with
$a_1+\ldots+ a_n=k$ and $0\le a_j\le s$, $1\le j\le n$, the
number of the corresponding elementary boxes
$\Dl^M_A\in\gE_s(q,n)$ with $M=(m_1,\ldots,m_n)\in\BN^n_0$,
$m_j\in\{0,1,\ldots, q^{a_j}-1\}$, $1\le j\le n$, is equal to
$q^k$, that coinsides with the cardinality of $D$. Therefore,
each of these elementary boxes $\Dl^M_A\in\gE_s(q,n)$ of volume
$q^{-k}$ contains exactly one point of $D$. Thus, $D$ is an
optimum $[ns,k]_s$-distribution. $\square$
\medskip

{\bf 2.5. Additional Remarks.} We complete the present section
with few remarks.

\medskip
{\it Remark 2.1.} In the special of $s=1$ the metrics
$\rho$ and $\ka$ coinside. Hence, the results and methods of the
present paper contain those for the Hamming metric $\ka$ as a
specific case. In particular, for $s=1$ Definition~2.1 coincides
with the known definition of MDS codes in 
the Hamming metric (cf. [16, Chapter 11]). 

It should be pointed out that both metrics $\varkappa$ and
$\rho$ belong to a broad class of metrics of the following kind.
Let $E$ be a vectors space of finitely many dimensions over an
arbitrary field. Fix a collection $\Cal M$ of subspaces $V\s E$
such that $E\in\Cal M$, $\{0\}\in\Cal M$, and for any two subspaces
$V_1,V_2\in\Cal M$ their vector sum $V_1+V_2\in\Cal M$. Now we put 
$$
\si(\Cal M;X)=\min\{\dim V:X\in V\in \Cal M\},\quad X\in E. 
\tag2.26
$$

It is obvious that $\si(\Cal M;\cdot)$ is a metric on $E$. 

Certainly, the metric (2.26) may be trivial for a bad choice of
the collection $\Cal M$. For example, $\si(\Cal M;X)=1$ for all $X\in
E\setminus\{0\}$ if $\Cal M$ contains each one dimensional subspace
of $E$. 

However, specific choises of $\Cal M$ lead to very interesting
metrics. For example, the metric $\rho$ can be defined by (2.26)
for $E=\Mat_{n,s}(\BF_q)$ and $\Cal M$ consisting of subspaces
$V^0_A\s \Mat_{n,s}(\BF_q)$ given by Lemma~2.2. 

A rich theory of MDS codes in metrics (2.26) can be developed
for approproate choises of collections $\Cal M$. Furthermore,
elements of such MDS codes turn out to be uniformly distributed
in a set of affine subspaces in $E$ parallel to subspaces $V\in
\Cal M$. Curiously enough, the indicated combinatorial
constructions can be given over arbitrary fields, including the
usual fields of real and complex numbers. 
\medskip

{\it Remark 2.2.} It should be of interest to determine
transformations preserving the metric $\rho$. Obviously,
the symmetric group $S_n$ of all  permutations of rows of a matrix
$\Om\in\Mat_{n,s}(\BF_q)$ preserves the weight $\rho(\Om)$ (and
the weight $\ka(\Om)$ as well). Examples of another kind can be
given as follows. 

Let $T_s$ denote a group of all nonsingular lower triangular $s$
by $s$ matrices over $\BF_q$ (with arbitrary nonzeroth diagonal
elements). Then, it is readily seen that $\rho(\om v)=\rho(\om)$
for $\om=(\xi_1,\ldots,\xi_s)\in\Mat_{1,s}(\BF_q)$ and $v\in
T_s$. Therefore, the product 
$$
H=S_n\times T_s\times\ldots \times T_s=S_n\times T^n_s
\tag2.27
$$
forms a group of linear transformations preserving the weight
$\rho$. 

We shall return to the group (2.27) once again in Sec.~4.4 below
in the contex of MacWilliams-type theorems for the metric
$\rho$. 

\medskip
{\it Remark 2.3.} It is instructive to compare the metric $\rho$
adoption in coding theory and the theory of uniform
distributions. It is well known that the major field
of application for coding theory is related with the
transmission of information (cf. [16]). Whether the metric
$\rho$ can be used in this area? Very schematically the problem
can be outlined as follows. 

Suppose that we wish to transmit a sequence
$\Om=(\eta_1,\eta_2,\ldots)$ of digits $\eta_j\in\BF_q$ across a
noisy channel. If the channel noise generates equiprobable
errors in all possible positions, the concept of the Hamming
metric is quite adequate to the problem. However, the Hamming
metric turns out to be very crude if possible errors form
patterns of a specific shape. 

Consider, for example, a channel where errors are generated by a
periodic  spike-wise perturbation of period $s$. We split the
sequence of digits $\Om$ of length $ns$ to $n$ blocks of length
$s$: $\Om=(\om_1,\ldots,\om_n)$, $\om_j\in\Mat_{1,s}(\BF_q)$,
$1\le j\le n$. In each block
$\om_j=(\xi^{(j)}_1,\ldots,\xi^{(j)}_s)$ the errors have a
tendency to occur at the left end, just after a spike puls,
whereas the probability of errors at the right end, like
$(0,\ldots,0,1)$, is negligible.

Now, it is obvious that the non-Hamming metric defined in (2.14)
depicts the actual structure of the indicated problem much more
adequately. Moreover, a code $\CC\subset\Mat_{n,s}(\BF_q)$ with
a large weight $\rho(\CC)$ is capable of detecting and
correcting all the 
channel errors $\Om'\in\Mat_{n,s}(\BF_q)$ with
$2\rho(\Om')<\rho(\CC)$, even though its Hamming weight
$\varkappa(\CC)$ is small (see (2.20)).

We refer to [25] for another possible application of the metric
$\rho$ to the transmission of information. 
Presumably, efficient decoding algorithms can be also given for
MDS codes in the metric $\rho$ (cf. [16, Sections~9.6 and
10.10]). Certainly, the study of such algorithms may be of
interest only if the metric $\rho$ is really practicable for
realistic systems of communication. 

The foregoing implies in passing that our basic parameters $n$
and $s$ play distinct roles is applications of the metric $\rho$
to coding theory and the theory of uniform distributions. In the
former case the parameter $s$  is given by constructional
features of the channel, while the parameter $n$ is free and can
be chosen arbitrary large if we wish to transmit long messages. 

Exactly the converse situation arises from problems related with
uniform distributions. Here $s\to\infty$ while $n$ is fixed,
because we are interested in more and more uniformly distributed
point sets in the unit cube of a given dimension $n$. 

A comprehensive study of uniform bounds for main code
characteristics at arbitrary values of the parameters $n$ and
$s$ should be of great importance. 

\medskip

\centerline{\bf 3. Weight Spectra}
\medskip

{\bf 3.1. Spheres and Balls in the Metric $\rho$.}

\medskip
{\it Definition 3.1.}  Given a code $\CC\subset
\Mat_{n,s}(\BF_q)$. The following set of $ns$ nonnegative
integers 
$$
w_r(\Om')=w_r(\CC;\Om')=\#
\{\Om\in\CC:\rho(\Om-\Om')=r\},\quad r\in\BN_0,\quad
0\le r\le ns,
\tag3.1
$$
is called the {\it weight spectrum} of the code $\CC$ relative
to an element $\Om'\in\CC$. Similarly, for a distribution
$D\subset \BQ^n(q^s)$ the set of nonnegative integers 
$$
\gathered
w_r(X')=w_r(D;X')=\#\{
X\in D:\rho(X\ominus X')=r\},
\\
r\in\BN_0,\quad 0\le r\le ns,
\endgathered
\tag3.2
$$
is called the {\it weight spectrum} of the distribution $D$
relative to a point $X'\in D$.
\medskip

It is obvious that weight spectra of distributions and the
corresponding codes coinside: $w_r(D;X')=w_r(\CC\l
D\r;\Om(X'))$. 

A remarkable thing is that weight spectra of all MDS codes in
the metric $\rho$ as well as optimum distributions can be evaluated
precisely. For MDS codes in the Hamming metric such results are
well known in coding theory (see [16, Chapter 11]).

Clearly,  weights (3.1) and (3.2) are equal to the number of
elements $\Om\in\CC$ or points $X\in D$ lying on spheres in the
metric $\rho$. Hence, we need to consider these geometric
subjects closer. Here we consider only the space $\BQ^n(q^s)$.
Certainly, the corresponding definitions and results for the
space $\Mat_{n,s}(\BF_q)$ can easily be given as well. 
Let 
$$
\goth S(r)=\{X\in\BQ^n(q^s):\rho(X)=r\},\quad
0\le r\le ns,
\tag3.3
$$
be a sphere of radius $r\in\BN_0$ in the metric $\rho$ on the
space $\BQ^n(q^s)$.

Split the unit interval into union of the following $s+1$
subintervals 
$$
[0,1)=\bigcup_{0\le b\le s}g_b,
\tag3.4
$$
where $g_0=[0,q^{-s})$, $g_1=[q^{-s},q^{-s+1}),\ldots$,
$g_s=[q^{-1},1]$. Then, the unit cube splits into a union of the
following $(s+1)^n$ disjoint boxes
$$
U^n=\bigcup_B G_B,\quad
G_B=\prod^n_{j=1} g_{b_j},
\tag3.5
$$
where $B=(b_1,\ldots,b_n)\in\BN^n_0$, $0\le b_j\le s$, $1\le
j\le n$. The boxes $G_B$  are said to be {\it fragments}. 

\medskip{\bf Lemma 3.1.} {\it Each sphere $(3.3)$ is a disjoint
union of the following fragments
$$
\gS(r)=\bigcup_{b_1+\ldots +b_n=r} G_B.
\tag3.6
$$
}
\medskip

{\it Proof.} Comparing  $q$-ary representation (2.5) with the
definitions of the $\rho$-weight (2.11) and subintervals $g_b$
in (3.4), we observe the following: a point $x\in \BQ(q^s)$
belongs to a subinterval $g_b$ if and only if its weight
$\rho(x)=b$. Thus, a point $X=(x_1,\ldots,x_n)^T\in\BQ^n(q^s)$
falls onto a sphere $\gS(r)$ if and only if $X$ belongs to one
of fragments $G_B$ with $b_1+\ldots+b_n=r$. $\square$
\medskip

For a given integer vector $B=(b_1,\ldots,b_n)\in\BN^n_0$ we
denote by $\ka(B)$ the number of nonzeroth coordinates $b_j$,
$1\le j\le n$, thus, $\ka(B)$ is equal to the Hamming weight of
the vector $B$. Consider a subset $\gS_l(r)\subset \gS(r)$
consistings of fragments $G_B\subset \gS(r)$ with $\ka(B)=l$, 
$$
\gS_l(r)=\bigcup_{\Sb b_1+\ldots+b_n=r,\\
\ka(B)=l\endSb} G_B,\quad l\in\BN_0,\quad
0\le l\le n.
\tag3.7
$$

From (3.6) and (3.7) we obtain that
$$
\gS(r)=\bigcup_{0\le l\le n}\gS_l(r).
\tag3.8
$$
Certainly, some of subsets $\gS_l(r)$ in (3.8) may be empty. 

We introduce nonnegative integers $\si_s(l,r)$, $r\in\BN_0$,
$l\in\BN_0$, $0\le l\le n$, by setting
$$
\gathered
\si_s(l,r)=\#\{
A=(a_1,\ldots,a_l)\in\BN^l:a_1+\ldots+ a_l=r,
\\
0<a_j\le s,1\le j\le l\}
\endgathered
\tag3.9
$$
for $l\ge1$, for $l=0$ we put $\si_s(0,0)=1$ and $\si_s(0,r)=0$
for $r\ge1$. In particular, for $s=1$ we have 
$$
\si_1(l,r)=\dl_{l,r},
\tag3.10
$$
where $\dl_{l,r}$ is the Kronecker symbol.

Moreover,  from  definition (3.9),
  we can easily  obtain the following asymptotics 
$$
\si_s(l,st+c)=\sqrt lv(l,t)s^{l-1}+O(s^{l-2})
\tag3.11
$$
as $s\to\infty$. In (3.11) $1\le t\le l-1$, $c\in\BN$ are
integers, and $v(l,t)$ is the $(l-1)$-dimensional volume of a
section of the $l$-dimensional unit cube $U^l$ by the hyperplane
$x_1+\ldots+x_l=t$. It is evident that $v(l,t)>0$. 

We write $\left(\matrix n\\
i\endmatrix\right)$ for the usual binomial coeffcients,
moreover, $\left(\matrix n\\
i\endmatrix\right)=0$ for $i>n$. 

\medskip{\bf Lemma 3.2.} {\it For each subset $\gS_l(r)$ the number of
fragments $G_B$ in  union $(3.7)$ is equal to $\left(\matrix n\\
l\endmatrix\right)\si_s(l,r)$. 
}
\medskip

{\it Proof.} Let $l=0$. Then, there is only one fragment
$G_B=G_0=[0,q^{-s})^n$ with $\ka(B)=0$. It is obvious that the
fragment $G_0$ is present in  union (3.7) only for $t=0$.
This proves Lemma~3.2 for $l=0$. 

Let $l\ge1$. Consider one of $\left(\matrix n\\
l\endmatrix\right)$ possible choices of $l$ indices
$J=\{j_1,\ldots, j_l\}\subset\{1,\ldots,n\}$. It is obvious (cf.~
Lemma~3.1) that the numbers of fragments $G_B\subset\gS_l(r)$
with such $B$ that $0<b_j\le s$ for $j\in J$ and $b_j=0$ for
$j\not\in J$ is equal to $\si_s(l,t)$. This proves Lemma 3.2 for
$l\ge1$. $\square$
\medskip

Using Lemma~3.2 and  representations (3.7) and (3.8), we can
easily evaluate the number of points $X\in\BQ^n(q^s)$ lying on a
sphere $\gS(r)$. First of all, we notice that 
$$
\#\{\BQ(q^s)\cap g_i\}=
\cases 1&\quad\text{if}\quad i=0,\\
(q-1)q^{i-1}&\quad\text{if}\quad i=1,\ldots,s.
\endcases
$$
Now, from (3.5) we conclude that
$$
\#\{\BQ^n(q^s)\cap G_B\}=
(q-1)^{\ka(B)}q^{b_1+\ldots+b_n-\ka(B)}.
$$
Therefore, using Lemma~3.2 and (3.7), (3.8), we obtain the
relation 
$$
\#\{\BQ^n(q^s)\cap\gS(r)\}=
\sum^n_{l=0} \left(\matrix n\\
l\endmatrix\right)\si_s(l,r)(q-1)^l q^{r-l}.
\tag3.12
$$
Notice that the terms with $l>r$ are absent in (3.12), because
$\si_s(l,r)=0$ if $l>r$ (cf. (3.9)).
Let 
$$
\gB(t)=\{X\in\BQ^n(q^s):\rho(X)\le t\}
\tag3.13
$$
be a ball of radius $t\in\BN_0$ in the metric $\rho$. It is
obvious that 
$$
\gB(t)=\bigcup_{0\le r\le t} \gS(r).
\tag3.14
$$

Let $D\subset\BQ^n(q^s)$ be  a distribution of $N$ points.
Assume that its weight $\rho(D)\ge 2t+1$, $t\in\BN_0$ (see
(2.20)). Then the balls $\{\gB(t)\oplus X,X\in D\}$ are pairwise
disjoint. This observation together with relations (3.12),
(3.13), and (3.14) implies the following {\it ball packing
bound} (cf. [25, Theorem~3])
$$
N\sum^t_{r=0}\sum^n_{l=0}\left(\matrix n\\
l\endmatrix\right) \si_s(l,t)(q-1)^l q^{t-l}\le
q^{ns}.
\tag3.15
$$

In view of (3.10) relation (3.15) for $s=1$ coinsides with
the well know Hamming bound for usual codes (cf. [9, 15, 16]).
In a similar spirit known results of coding theory (like the
Gilbert--Varshamov and Plotkin bounds) can be also given with
respect to the 
metric $\rho$ (cf. [25, Theorems~2 and 4]).
 In connection with this, we should mention a very
interesting combinatorial problem on the description of {\it
perfect distributions} $D\subset \BQ^n(q^s)$ (and the
corresponding {\it perfect codes} $\CC\l D\r\subset
\Mat_{n,s}(\BF_q)$ in metric $\rho$) which provide exact
close-packings (without gaps and overlappings) of the space
$\BQ^n(q^s)$ by balls $\{\gB(t)\oplus X,X\in D\}$ or spheres
$\{\gS(r)\oplus X,X\in D\}$. However, in the present paper we
shall not consider these questions closer. 

The following result gives a description of balls (3.13) in
terms of elementary boxes (1.5). 

\medskip{\bf Lemma 3.3.} {\it Each ball $(3.13)$ is a union of the
following elementary boxes $(1.5)$ from the collection $(1.8)$ 
$$
\gB(t)=\bigcup_{a_1+\ldots+a_n=t}\Dl^0_{A^*},
\tag3.16
$$
where $A^*=(a^*_1,\ldots,a^*_n)$ with $a^*_j=s-a_j$, $0\le
a_j\le s$, $1\le j\le n$.
}
\medskip

{\it Proof.} Our arguments here are quite similar to those in
the proof of Lemma~2.3. Let $X\in \gB(t)$. Writing 
inequality (2.24) for all coordinates of the point $X$, we
conclude that $X\in\Dl^0_{A(X)}$ with
$A(X)=(a_1(x_1),\ldots,a_n(x_n))$ where $a_j(x_j)=s-\rho(x_j)$,
$1\le j\le n$, and $\rho(x_1)+\ldots+\rho(x_n)=\rho(X)\le t$.
Hence, there exists an elementary box $\Dl^0_A$ in the union
(3.16) such that $\Dl^0_A\supseteq\Dl^0_{A(X)}\ni X$. 

Let $X\not\in\gB(t)$, that is, $\rho(X)>t$. Writing 
inequalities (2.13) and (2.25) correspondingly for zeroth and
nonzeroth coordinates of the point $X$, we observe the
following: if the point $X$ falls into an elementary box
$\Dl^0_{A^*}\in \gE_s(q,n)$ with $A^*=(a^*_1,\ldots,a^*_n)$ and
$a^*_j=s-a_j$, $0\le a_j\le s$, $1\le j\le n$, then
$\rho(x_j)\le a_j\le s$, $1\le j\le n$. Therefore, we have
$a_1+\ldots+a_n\ge\rho(X)>t$, and hence, the elementary box
$\Dl^0_{A^*}$ does not belong to  union (3.16). $\square$
\medskip

{\bf 3.2. Weight Spectra of MDS Codes and Optimum
Distributions.} 
Our main result on weight spectra can be stated as follows. 

\medskip{\bf Theorem 3.1.} {\it Let $D\s\BQ^n(q^s)$ be an optimum
$[ns,k]_s$-distribution and $\CC=C\l D\r$ the corresponding MDS
$[ns,k]_s$-code in the metric $\rho$. Then weights $(3.1)$ and
$(3.2)$ are independent of elements $\Om'\in\CC$ and points
$X'\in D$, and $w_r=w_r(\Om')=w_r(X')$ are given by
$$
w_0=1,\quad w_r=0
\tag3.17
$$
for $0\le r<\rho(\CC)= ns-k+1$, and
$$
w_r=\sum^n_{l=1}\left(\matrix n\\
l\endmatrix\right)\si_s(l,r)\sum^{r-\rho(\CC)}_{t=0}
(-1)^t \left(\matrix l\\
t\endmatrix\right)(q^{r-\rho(\CC)+1-t}-1)
$$
$$
=(q-1)\sum^n_{l=1}\left(\matrix n\\
l\endmatrix\right)\si_s(l,r)\sum^{r-\rho(\CC)}_{t=0}
(-1)^t \left(\matrix l-1\\
t\endmatrix\right)q^{r-\rho(\CC)-t}
\tag3.18
$$
for $\rho(\CC)\le t\le ns$. 
}
\medskip

We wish to discuss some corollaries of Theorem~3.1. First of
all, we notice that in view of (3.10) relations (3.18) for
$s=1$ imply 
$$
w_r=\sum^{r-\rho(\CC)}_{t=0}(-1)^t\left(\matrix l\\
t\endmatrix\right)(q^{r-\rho(\CC)+1-t}-1)
$$
$$
=(q-1)\sum^{r-\rho(\CC)}_{t=0}(-1)^t\left(\matrix l-1\\
t\endmatrix\right)q^{r-\rho(\CC)-t}.
$$
These relations coinside with known identities (cf. [16, Chapter
11, Theorem 6]) for weight spectra of MDS codes in the Hamming
metric (see Remark 2.1). 

Now we consider the special case of $k=s$. Introduce the
following positive integers
$$
\gathered
\tsi_s(n,r)=\#\{
A=(a_1,\ldots,a_n)\in\BN^n_0:
a_1+\ldots+a_n=r,
\\
0\le a_j\le s,\,1\le j\le n\}.
\endgathered
\tag3.19
$$
Comparing the definitions (3.9) and (3.19), we immediately
obtain the relation
$$
\tsi_s(n,r)=\sum^n_{l=1}\left(\matrix n\\
l\endmatrix\right)\si_s(l,r),\quad
r\ge1.
\tag3.20
$$

Theorem 3.1 implies the following result. 

\medskip{\bf Theorem 3.2.} {\it Let $D$ be an optimum
$[ns,s]_s$-distribution or, that is the same (cf.
Proposition~$1.1$), a $(0,s,n)$-net. Then weights $(3.2)$ are
given by 
$$
w_0=1,\quad w_r=0
\tag3.21
$$
for $0\le r<\rho(D)=(n-1)s+1$, and
$$
w_r=\tsi_s(n,r)
\sum^{r-\rho(D)}_{t=0}(-1)^t\left(\matrix n\\
t\endmatrix\right)(q^{r-\rho(D)+1-t}-1)
$$
$$
=\tsi_s(n,r)(q-1)\sum^{r-\rho(D)}_{t=0}(1-1)^t\left(\matrix n-1\\
t\endmatrix\right) q^{r-\rho(D)-t}
\tag3.22
$$
for $\rho(D)\le r\le ns$. In particular,
$$
w_r=\tsi_s(n,r)(q-1)^nq^{r-\rho(D)-n+1}
\tag3.23
$$
for $\rho(D)+n-1\le r\le ns$.
}
\medskip

{\it Proof.} Evidently,  relations (3.17) imply (3.21). Let
$\rho(D)=(n-1)s+1\le r\le ns$, then $\si_s(l,r)=0$ for $1\le
l\le n-1$ and $\si_s(n,r)=\tsi_s(n,r)$ (see (3.9), (3.19), and
(3.20)). Hence,  relations (3.18) imply (3.22). Finally, 
relations (3.23) follow at once from (3.22), because, for
$\rho(D)+n-1\le r\le ns$ the summation in (3.22) is taken over
all indices $0\le t\le n-1$. $\square$
\medskip

By definition, weights (3.1) and (3.2) are nonnegative. Together
with Theorems~3.1 and 3.2 this leads to a series of necessary
conditions for the existence of optimum distributions and MDS
codes. From  relation (3.17) we find the first nontrivial terms
in the weight spectrum
$$
w_{\rho(D)}=\tsi_s(n,\rho(D))(q-1),
\tag3.24
$$
and
$$
w_{\rho(D)+1}=(q-1)\sum^n_{l=1}\left(\matrix n\\
l\endmatrix\right)(q+1-l)\si_s(l,\rho(D)+1).
\tag3.25
$$
Notice that  relation (3.20) has been used in (3.24).
Moreover, writing  relation (3.25), we assume that
$\rho(D)+1\le ns$; this means that $s>1$ if $k=s$. 

For $k=s$  relation (3.25) takes the form 
$$
w_{\rho(D)+1}=(q-1)\si_s(n,\rho(D)+1)(q-n+1),
\tag3.26
$$
because $\si_s(l,(n-1)s+2)=0$ for $l<n$ (see (3.9)).

Since the weight (3.26) is nonnegative, we obtain the following
necessary condition for the existence of optimum
$[ns,s]_s$-distributions:  
$$
q\ge n-1.
\tag3.27
$$

For $(0,n,s)$-nets  condition (3.27) is well known (cf. [4,
Lemma~3.28] and [20, Corollary 5.11]). The above arguments can
be considered as an interpretation of  condition (3.27) in
terms of weight spectra. 

From relation (3.25) we can also derive the following
general result. 

\medskip{\bf Proposition 3.1.} {\it Suppose that $k=st$ for given $1\le
t\le 
n-1$, and for an infinite sequence of $s\in\BN$ there exist optimum
$[ns,k]_s$-distributions (or, that is the same, MDS
$[ns,k]_s$-codes in the metric $\rho$). Then  condition
$(3.27)$ holds. 
}
\medskip

{\it Proof.} Substituting  asymptotics (3.11) to relation
(3.25), we find that 
$$
(q-1)(q+1-n)\sqrt n \,v(n,n-t)s^{n-1}
+O(s^{n-2})\ge0,\quad
s\to\infty,
$$
with the constant $v(n,n-t)>0$, and  inequality (3.27)
follows. $\square$
\medskip

\medskip
{\it Remark 3.1.} For given $1\le k\le ns$, $q\ge n-1$
and arbitrary $s\in\BN$ optimum $[ns,k]_s$-distributions will be
explicitly constructed in Section~5 (see Theorem~5.2 below).
Hence,  condition (3.27) is actually necessary and sufficient
for the existence of infinite sequences of optimum distributions
and MDS codes in the metric $\rho$.

\medskip
{\bf 3.3. Proof of Theorem 3.1.}
First of all, we notice that Theorem 3.1 may be proved only for
optimum distributions, because a given MDS $[ns,k]_s$-code
$\CC\s\Mat_{n,s}(\BF_q)$ can be identified with the
corresponding optimum $[ns,k]_s$-distribution $D=D\l
\CC\r\s\BQ^n(q^s)$ by Theorem~2.1. 

Next, from Definition 2.1 and Theorem~2.1 we conclude that if
$D\s\BQ^n(q^s)$ is an optimum $[ns,k]_s$-distribution, then so is
any distribution $D\oplus X'$ shifted by an arbitrary vector
$X'\in\BQ^n(q^s)$. Thus, without loss of generality, in proving
of Theorem~3.1 we may assume that $0\in D$ and consider 
weight spectrum (3.2) only relative to the origin $X'=0$. In
so doing, we put $w_r=w_r(0)$. 

Notice also that  relations (3.17) hold trivially, in view of
Definition~2.1 and Theorem~2.1. Thus, we have to consider only
weights $w_r$ with $r\ge\rho(D)= ns-k+1$. 

With the above conventions we turn to the proof of  relations
(3.18). From Lemma~3.1 we obtain the relation 
$$
w_r=\#\{
D\cap\gS(r)\}=
\sum_{b_1+\ldots+b_n=r}\#\{D\cap G_B\}.
\tag3.28
$$

Consider a fragment $G_B$, $B=(b_1,\ldots,b_n)$, in (3.28) with
$\ka(B)=l\ge1$. Notice that a fragment $G_B$ with $\ka(B)=0$ is
absent in (3.28) since $r>0$. Suppose that $b_j\ne 0$ $(0<b_j\le
s)$ for $j\in J$ and $b_j=0$ for $j\not\in J$, where
$J=\{j_1,\ldots,j_l\}
\subseteq\{1,\ldots,n\}$ is a subset of $l$ indices. Using
(3.4), (3.5), we conclude that the fragment $G_B$ has the
following representations in terms of elementary boxes (1.5) 
$$
G_B=\Dl^0_{A_0}\setminus \bigcup_{i\in J}\Dl^0_{A_i},
\tag3.29
$$
where $A_0=(a^{(0)}_1,\ldots,a^{(0)}_n)$ with 
$$
a_j=s-b_j,\quad 1\le j\le n,
\tag3.30
$$
and $A_i=(a^{(i)}_1,\ldots,a^{(i)}_n)$ with
$$
a^{(i)}_j=
\cases
s\quad&\text{if}\quad j\not\in J,\\
s-b_j\quad&\text{if}\quad j\in J\quad \text{and}\quad j\ne i,\\
s-b_i+1\quad&\text{if}\quad j\in J\quad\text{and}\quad j=i.
\endcases
\tag3.31
$$
Obviously, we have the equality
$$
a^{(0)}_1+\ldots+a^{(0)}_n=ns-b_1-\ldots-b_n=ns-r.
\tag3.32
$$

Let $I=\{i_1,\ldots,i_t\}\subseteq J=\{ j_1,\ldots,j_l\}$ be a
subset of $t$ indices. Then, 
$$
\Dl^0_{A_{i_1}}\cap\ldots\cap \Dl^0_{A_{i_t}}=\Dl^0_{A_I},
\tag3.33
$$
where $A_I=(a^{(I)}_1,\ldots,a^{(I)}_n)$ with 
$$
a^{(I)}_j=\max \{a^{(i_1)}_j,\ldots,a^{(i_t)}_j\}
=
\cases
s-b_j\quad&\text{if}\quad j\not\in I,\\
s-b_j+1\quad&\text{if}\quad j\in I.
\endcases
\tag3.34
$$
Notice that  relations (3.31) have been used in (3.34).
Obviously, we have the equality 
$$
a^{(I)}_1+\ldots+a^{(I)}_n=ns-r+l.
\tag3.35
$$

Applying the principle of inclusion and exclusion to the summand
$\#\{D\cap G_B\}$ in (3.28), and using  relations (3.29),
(3.32), we obtain the formula 
$$
\#\{D\cap G_B\}=\{D\cap \Dl^0_{A_0}\}-
\sum^l_{t=1}\sum_{\Sb I\subset J,\\
\#\{I\}=t\endSb}(-1)^t\#
\{D\cap\Dl^0_{A_I}\},
\tag3.36
$$
where the inner sum is taken over all $\left(\matrix l\\
t\endmatrix\right)$ subsets $I=\{i_1,\ldots,i_t\}\subseteq
J=\{j_1,\ldots,j_l\}$ of $t$ indices. 

From  relations (3.30) and (3.34) we conclude that all
elementary boxes $\Dl^0_{A_0}$ and $\Dl^0_{A_I}$ in (3.36)
belong to  collection (1.8). Now we observe the following.
Let an elementary box $\Dl^0_A\in\gE_s(q,n)$, then 
$$
\#\{D\cap\Dl^0_A\}=
\cases
q^{s-a_1-\ldots-a_n}\quad&\text{if}\quad a_1+\ldots+a_n<k,\\
1\quad&\text{if}\quad a_1+\ldots+a_n\ge k.
\endcases
\tag3.37
$$
The first equality in (3.37) follows from Lemma~1.2 (i), and the
second follows from Lemma~1.2 (ii) and the above assumption
that $0\in D$. 

Substituting (3.37) to (3.36), and using (3.32) and (3.35), we
find that 
$$
\#\{D\cap G_B\}=\sum^{r- ns-k-1}_{t=0}(-1)^t
\left(\matrix l\\
t\endmatrix\right) q^{r- ns-k-t}+\sum^l_{t=r- ns-k}
(-1)^t \left(\matrix l\\
t\endmatrix\right)
$$
$$
=\sum^{r-\rho(D)}_{t=0}(-1)^T \left(\matrix l\\
t\endmatrix\right)
(q^{r-\rho(D)+1-t}-1)
$$
$$
=
(q-1)\sum^{r-\rho(D)}_{t=0}(-1)^t \left(\matrix
l-1\\
t\endmatrix\right)q^{r-\rho(D)-t}.
\tag3.38
$$
Notice that the well know identities for binomial coefficients:
$$
\sum\limits^l_{t=0}(-1)^t\left(\matrix l\\
t\endmatrix\right)=0\quad \text{and}\quad
\left(\matrix l\\
t\endmatrix\right)=\left(\matrix l-1\\
t\endmatrix\right)+\left(\matrix l-1\\
t-1\endmatrix\right)
$$
 have been used in  calculation (3.38). 

Relation (3.28) can be written  in the form 
$$
w_r=\sum^n_{l=1}\sum_{B:\ka(B)=l}\#\{D\cap G_B\},
\tag3.39
$$
where the inner sum is taken over all fragments $G_B\s \gS(r)$
with $\ka(B)=l$. The number of such fragments is equal to $\left(\matrix l\\
t\endmatrix\right)\si_s(l,r)$ by Lemma~3.2. 

Substituting (3.38) to (3.39), we obtain the relation in
question (3.28).

The proof of Theorem~3.1 is complete. $\square$

\medskip
\centerline{\bf 4. Linear Codes and Distributions}
\medskip

{\bf 4.1. Duality for Codes and Distributions.}

\medskip{\it Definition 4.1.}  A subset
$\CC\s\Mat_{n,s}(\BF_q)$ is 
called a {\it linear code} if $\CC$ is a subspace in
$\Mat_{n,s}(\BF_q)$. Similarly, a subset $D\s\BQ^n(q^s)$ is
called a {\it linear distribution} if $D$ is a subspace in the
vector space $\BQ^n(q^s)$ over $\BF_q$ with respect to the
arithmetic operations $(2.10)$.
\medskip

For linear codes and distributions the definitions of their
weights (2.17) and (2.18) take the form 
$$
\la(\CC)=\min\{\la(\Om):\Om\in\CC\setminus\{0\}\},\quad
\la(D)=\min\{\la(X):X\in D\setminus\{0\}\},
\tag4.1
$$
where $\la$ is one of the metrics $\ka$ or $\rho$. 

Notice that  relations (4.1) are quite similar to the known
definition of the homogeneous minimum of a lattice in Euclidean
spaces (cf. [9] and [16]). 

We introduce the following {\it inner product}
$\l\boldsymbol\cdot,\boldsymbol\cdot\r$ on the space
$\Mat_{n,s}(\BF_q)$. At first, let $n=1$ and
$\om_1=(\xi'_1,\ldots,\xi'_s)$,
$\om_2=(\xi''_1,\ldots,\xi''_s)\in\Mat_{1,s}(\BF_q)$, then we
put 
$$
\l \om_1,\om_2\r=
\l\om_2,\om_1\r=\sum^s_{i=1}\xi'_i\xi''_{s+1-i}.
\tag4.2
$$
Now, let
$\Om_i=(\om^{(1)}_i,\ldots,\om^{(n)}_i)^T\in\Mat_{n,s}(\BF_q)$,
$i=1,2$, where $\om^{(j)}_i\in\Mat_{1,s}(\BF_q)$, $1\le j\le n$.
Then we put 
$$
\l\Om_1,\Om_2\r=\l\Om_2,\Om_1\r=
\sum^n_{j=1}\l\om^{(j)}_1,\om^{(j)}_2\r.
\tag4.3
$$

Certainly, instead of (4.2) we could write the inner product in
a common form
$(\om_1,\om_2)=\xi'_1\xi''_1+\ldots+\xi'_s\xi''_s$. In doing so,
$(\om_1,\om_2)=\l\om_1,\Phi_0\om_2\r$ with a nonsingular $s$ by
$s$ matrix $\Phi_0$ over $\BF_q$. Recall that an arbitrary
non-singular bilinear form
$\Phi(\om_1,\om_2)=\l\om_1,\Phi\om_2\r$, $\det\Phi\ne0$, could
be used as an appropriate inner product on the space
$\Mat_{1,n}(\BF_q)\simeq \BF^n_q$ (cf. [14, Chapters 13, 14]).
Definition (4.2) is preferable here, because otherwise, the
matrix $\Phi_0$ would be involved in all subsequent
considerations. 

Using  relation (2.8), we can transfer  inner product
(4.3) onto the vector space of rational $q$-ary points
$\BQ^n(q^s)$ by setting
$$
\l X_1,X_2\r=\l\Om(X_1),\Om(X_2)\r,\quad X_1,X_2\in\BQ^n(q^s). 
\tag4.4
$$

\medskip{\it Definition 4.2.}  Given a linear code
$\CC\s\Mat_{n,s}(\BF_q)$,  the {\it dual code}
$\CC^\perp\s$\break $\s\Mat_{n,s}(\BF_q)$ is defined by 
$$
\CC^\perp=\{\Om_2\in\Mat_{n,s}(\BF_q):\l\Om_1,\Om_2\r=0\quad\text{for
all}\quad \Om_1\in\CC\}.
\tag4.5
$$
Similarly, for a given linear distribution $D\s\BQ^n(q^s)$ the
{\it dual distribution} $D^\perp\s\BQ^n(q^s)$ is defined by
$$
D^\perp=\{X_2\in\BQ^n(q^s):\l X_1,X_2\r=0\quad\text{for all}\quad
X_1\in D\}.
\tag4.6
$$
\medskip

The following simple corollaries of Definition~4.2 should be
mentioned. Dual codes (4.5) and distributions (4.6) are also
linear subspaces. Moreover, 
$$
(\CC^\perp)^\perp=\CC,\quad (D^\perp)^\perp=D,
\tag4.7
$$
and (cf. (2.9))
$$
\CC\l D^\perp\r=(\CC\l D\r)^\perp,\quad
D\l \CC^\perp\r=(D\l \CC\r)^\perp.
\tag4.8
$$

Thus, in view of (4.7), linear codes
$C,C^\perp\s\Mat_{n,s}(\BF_q)$ and distributions
$D,D^\perp\s\BQ^n(q^s)$ are pairs of mutually dual subspaces. 

Let $\CC\s\Mat_{n,s}(\BF_q)$ and $D\s\BQ^n(q^s)$ be subspaces of
dimension $k$, then $\CC^\perp$ and $D^\perp$ are subspaces of
dimention $k^\perp=ns-k$. This implies in passing the following
relations for their cardinalities
$$
\#\{\CC\}=\#\{D\}=q^k,\quad
\#\{\CC^\perp\}=\#\{D^\perp\}=q^{k^\perp},
\tag4.9
$$
and
$$
\#\{\CC\}\#\{\CC^\perp\}=\#\{D\}\#\{D^\perp\}=q^{ns}.
\tag4.10
$$
Furthermore, let $\CC$ and $D$ be linear an $[ns,k]_s$-code and,
correspondingly, an $[ns,k]_s$-distribution, then $\CC^\perp$ and
$D^\perp$ are linear an $[ns,ns-k]_s$-code and, correspondingly,
$[ns,ns-k]_s$-distribution. 

Our results on metric properties of mutually dual codes and
distributions can be stated as follows. 

\medskip{\bf Theorem 4.1.} {\rm(i)} {\it  Let
$\CC\s\Mat_{n,s}(\BF_q)$ be a linear MDS $[ns,k]_s$-code in the
metric $\rho$. Then 
$\CC^\perp$ is a linear MDS $[ns,ns-k]_s$-code in the metric
$\rho$. 

{\rm(ii)} Let $D\s\BQ^n(q^s)$ be a linear optimum
$[ns,k]_s$-distribution. Then $D^\perp$ is a linear optimum
$[ns,ns-k]_s$-distribution. 
}
\medskip

It is notable that linear $(\dl,n,s)$-nets (cf. Definition~1.1)
can also be characterized completely in terms of their dual
distributions (cf. [20, Theorems~6.10 and 6.14]).  

\medskip{\bf Theorem 4.2.} {\it Let $D$, $D^\perp\s\BQ^n (q^s)$ be
mutually dual {\it linear} distributions of dimensions $k=s$ and
$k^\perp=(n-1)s$, respectively. Then the following two statements
are equivalent: 

{\rm(i)} $D$ is a $(\dl,n,s)$-net of deficiency $\dl$, 

{\rm (ii)} $D^\perp$ has a weight
$\rho(D^\perp)\ge s+1-\dl$.
}
\medskip

It should be noted that the statements of Theorem~4.1 for $k=s$
and Theorem~4.2 for $\dl=0$ actually coincide (cf.
Proposition~1.1 and Theorem~2.1). 

For the proof of Theorems~4.1 and 4.2 we need to recall some
known facts about the Fourier transform on vector spaces over
finite fields. For details we refer to [15] and [16]. Notice
also that the necessary facts are given here in the form adapted
to the vector spaces $\BQ^n(q^s)$. 

\medskip
{\bf 4.2. Fourier Transform on the Space of Distributions.} 
Introduce a function $\Psi(X_1,X_2)$, $X_1,X_2\in\BQ^n(q^s)$, by
setting
$$
\Psi(X_1,X_2)=\exp\left(2\pi\sqrt{-1}\frac1p\,\Tr\l X_1,X_2\r\right),
\tag4.11
$$
where we write $\Tr\xi=\xi+\xi^p+\xi^{p^2}+\ldots+\xi^{p^{e-1}}$
for the trace of an element $\xi\in\BF_q$, $q=p^e$, over
$\BF_p$. 

It is obvious that  function (4.11) satisfies the relation 
$$
\Psi(X,X_1\oplus X_2)=\Psi(X,X_1)\Psi(X,X_2).
\tag4.12
$$
Hence, for a given $X$ the function $\Psi(X,\boldsymbol\cdot)$
is an additive character on the vector space $\BQ^n(q^s)$. One
can show (cf. [15, Section 5.1]) that each character of the
additive group of $\BQ^n(q^s)$ coincides with
$\Psi(X,\boldsymbol\cdot)$ for a suitable $X$. 

The Fourier transform $\tf$ of a complex-valued function
$f:\BQ^n(q^s)\to\BC$ is defined by 
$$
\tf(Y)=\sum_{X\in\BQ^n(q^s)}\Psi(Y,X)f(X).
\tag4.13
$$

The following results are well known in the theory of Abelian groups
(see, for example, [15, Chapters 5, 9] and [16, Chapter 5]). 

\medskip{\bf Lemma 4.1.} {\it Let $D$, $D^\perp\s\BQ^n(q^s)$ be
mutually dual linear distributions. Then,

{\rm (i)} One has the relation
$$
\sum_{X\in D}\Psi(Y,X)=\cases
\#\{D\}\quad&\text{if}\quad Y\in D^\perp,\\
0\quad&\text{if}\quad Y\not\in D^\perp.
\endcases
\tag4.14
$$

{\rm(ii)} One has the Poisson summation formula
$$
\sum_{X\in D}f(X)=q^{-ns}\#\{D\}\sum_{Y\in D^\perp}\tf(Y)
\tag4.15
$$
for an arbitrary function $f$.
}

Consider affine subspaces $V^M_A=\BQ^n(q^s)\cap\Dl^M_A$,
$\Dl^M_A\in\gE_s(q,n)$ as distributions in $\BQ^n(q^s)$ (see
Lemma~2.2). We write $\chi^M_A(\boldkey.)$ for the indicator
function of $V^M_A$, that is, 
$$
\chi^M_A(X)=\cases
1\quad&\text{if}\quad X\in V^M_A,\\
0\quad&\text{if}\quad X\not\in V^M_A.
\endcases
$$
Since $V^M_A=V^0_A\oplus Y^M_A$ with
$Y^M_A=\left(\frac{m_1}{q^{a_1}},\ldots,\frac{m_n}{q^{a_n}}\right)^
T\in\BQ^n(q^s)$, we have 
$$
\chi^M_A(X)=\chi^0_A(X\ominus Y^M_A).
\tag4.16
$$

We wish to find out distributions $(V^0_A)^\perp$ dual to
subspaces $V^0_A$ and evaluate the Fourier transform
$\tilde\chi^M_A(\cdot)$ of $\chi^M_A(\cdot)$. 

Notice that in view of (4.12), (4.13), and (4.16) we have the
relation 
$$
\tilde\chi^M_A(Y)=\Psi(Y,Y^M_A)\tilde\chi^0_A(Y).
\tag4.17
$$

\medskip{\bf Lemma 4.2.} {\it For a given
$A=(a_1,\ldots,a_n)\in\BN^n_0$, $0\le a_j\le s$, define the vector
$A^*=(a^*_1,\ldots,a^*_n)\in\BN^n_0$, $0\le a^*_j\le s$, $1\le
j\le n$,  by $a^*_j=s-a_j$, $1\le j\le n$. Then one
has the relations 
$$
(V^0_A)^\perp=V^0_{A^*},
\tag4.18
$$
and
$$
\tilde\chi^0_A(Y)=q^{ns-a_1-\ldots-a_n}\chi^0_{A^*}(Y).
\tag4.19
$$
}
\medskip

{\it Proof.} From Lemma~2.2 we derive that the subspace $V^0_A$
consists of points $X=(x_1,\ldots,x_n)^T\in\BQ^n(q^s)$ with
coordinates 
$$
x_j=\sum^{s-a_j}_{i=1}\xi_i(x_j)q^{i-s-1},\quad 1\le j\le n,
$$
where $\xi_i(x_j)$ are arbitrary elements of the field $\BF_q$.
Using  definitions (4.2), (4.3), (4.4), and (4.6), we
conclude that the dual subspace $(V^0_A)^\perp$ consists of
points $Y=(y_1,\ldots,y_n)^T\in\BQ^n(q^s)$ which satisfy the
equations 
$$
\xi_{s+1-i}(y_j)=0\quad\text{for}\quad 1\le i\le s-a_j,\quad
1\le j\le n,
$$
or in an equivalent form 
$$
\xi_i(y_j)=0\quad\text{for}\quad s\ge i\ge s+1-a^*_j,\quad
1\le j\le n.
$$
By Lemma~2.2 this implies  relation (4.18). 

Relation (4.19) follows at once from (4.18). It suffices to
write  relation (4.14) for $D=V^0_A$ and take into account
that $\#\{V^0_A\}=q^{ns-a_1-\ldots-a_n}$ (cf. relation (4.9)
and Lemma~2.2). $\square$
\medskip

\medskip{\bf Lemma 4.3.} {\it Let $D$, $D^\perp\s\BQ^n(q^s)$ be mutually
dual linear distributions. Then for each elementary box
$\Dl^M_A\in\gE_s(q,n)$ one has the relation
$$
\#\{D\cap\Dl^M_A\}=
q^{-a_1-\ldots-a_n}\#\{D\}\sum_{Y\in
D^\perp}\Psi(Y,Y^M_A)\chi^0_{A^*}(Y). 
\tag4.20
$$

In particular,
$$
\#\{D\cap\Dl^0_A\}=
q^{-a_1-\ldots-a_n}\#\{D\}\#\{D^\perp\cap \Dl^0_{A^*}\}.
\tag4.21
$$
}\medskip

{\it Proof.} It suffices to substitute (4.17) and (4.19) to the
Poisson summation formula (4.15) and use the relation 
$$
\#\{D\cap\Dl^M_A\}=\sum_{X\in D}\chi^M_A(X). \ \ \square
$$
\medskip

{\bf 4.3. Proof of Theorems 4.1 and 4.2.}
Theorems 4.1 and 4.2 are direct corollaries of the following
general result.

\medskip{\bf Proposition 4.1.} {\it Let $D$,
$D^\perp\s\BQ^n(q^s)$ be mutually dual linear distributions of
dimensions $d$ and 
$d^\perp=ns-d$, respectively. Then, for a given integer
$0\le\dl\le d$ the following two statements are equivalent:  

{\rm(i)} Each elementary box $\Dl^M_A\in\gE_s(q,n)$ of volume
$q^{-d+\dl}$ contains exactly $q^{\dl}$ points of $D$. 

{\rm(ii)} $D^\perp$ has a weight $\rho(D^\perp)\ge
ns-d^\perp-\dl+1=d-\dl+1$. 
}
\medskip

By  relation (4.8) the statement (i) and (ii) of Theorem 4.1
are equivalent. At the same time Theorem~4.1 (ii) follows at
once from Theorem~2.1 and  Proposition~4.1 for $d=k$ and
$\dl=0$. Similarly, 
Theorem~4.2 follows from Definition~1.1 and Proposition~4.1 for
$d=s$ and $\dl$ equal to the deficiency of the corresponding
net. 

The proof of Theorems~4.1 and 4.2 is complete. $\square$
\medskip

{\it Proof of Proposition 4.1.} (i) Let the statement (i) be
held. Then, writing  equality (4.21) for an elementary box
$\Dl^0_A\in\gE_s(q,n)$ of volume $q^{-d+\dl}$, and using 
relations (1.6) and (4.9), we conclude that $\#\{D^\perp\cap
\Dl^0_{A^*}\}=1$ for an arbitrary elementary box
$\Dl^0_{A^*}\in\gE_s(q,n)$ with
$a^*_1+\ldots+a^*_n=ns-d+\dl=d^\perp+\dl$. Hence each indicated
elementary boxes $\Dl^0_{A^*}$ contains a single point $Y=0$
belonging to $D^\perp$. By Lemma~3.3 this means that the ball
$\gB(t)\s\BQ^n(q^s)$ of radius $t=ns-d^\perp-\dl=d-\dl$ does not
contain points of $D^\perp$ other than the origin $Y=0$.
Therefore $\rho(D^\perp)\ge ns-d^\perp-\dl+1=d-\dl+1$ (cf. 
definition (4.1)). 

(ii) Let the statement (ii) be held. Then, by Lemma~3.3 each
elementary box $\Dl^0_{A^*}\in\gE_s(q,n)$ with
$a^*_1+\ldots+a^*_n=d^\perp+\dl=ns -d+\dl$ contains a single
point $Y=0$ belonging to $D^\perp$. Writing  equality (4.20)
for the indicated boxes $\Dl^0_{A^*}$, we find that
$\#\{D\cap\Dl^M_A\}= q^\dl$ for each elementary box
$\Dl^M_A\in\gE_s(q,n)$ of volume $q^{-d+\dl}$. $\square$
\medskip

{\bf 4.4. Enumerators of Mutually Dual Distributions.}
Relations (4.21) can be expressed in terms of the
corresponding generating functions. For a given distribution
$D\subset\BQ^n(q^s)$ we define its {\it box enumerator}
$\ph(D;z_1,\ldots,z_n)$, $(z_1,\ldots,z_n)\in\BC^n$, by 
$$
\ph(D;z_1,\ldots,z_n)=\sum_{\Sb 0\le a_j\le s,\\
1\le j\le n\endSb}\#\{D\cap\Dl^0_A\} z^{a_1}_1\ldots z^{a_n}_n 
$$
$$
=\sum_{X\in D}\,\sum_{\Sb 0\le a_j\le s,\\
1\le j\le n\endSb}\chi^0_A(X)z^{a_1}_1\ldots z^{a_n}_n. 
\tag4.22
$$
Expression (4.22) is a polynomial of degree $s$ in each
complex variable $z_j$. 

Substituting  relation (4.21) to (4.22), we immediately
obtain the following. 

\medskip{\bf Theorem 4.3.} {\it Let $D$ and $D^\perp\s\BQ^n(q^s)$ be
linear mutually dual distributions. Then, one has the relation 
$$
\ph(D;z_1,\ldots,z_n)
=
\#\{D\}\left(\frac{z_1\ldots z_n}{q^n}\right)^s\ph
\left(D^\perp;\frac q{z_1},\ldots,\frac q{z_n}\right).
\tag4.23
$$
}
\medskip

For a given linear distribution $D\s\BQ^n(q^s)$ we can also
introduce its {\it weight enumerator} $\Cal W(D;z)$, $z\in\BC$,
by 
$$
\Cal W(D;z)=\sum^{ns}_{r=0}w_r(D)\,z^r=\sum^{ns}_{r=0}
\#\{D\cap\goth S(r)\}z^r,
\tag4.24
$$
where $w_0(D)=1$, $w_1(D),\ldots,w_{ns}(D)$ is the weight
spectrum of $D$ in the metric $\rho$ (cf (4.1)), and $\goth
S(r)$ is  sphere (3.3). 

The Hamming metric weight enumerators of mutually dual codes
are related by the well known MacWilliams identity (see [16, Chapter
5, Theorem 13]). Is there an extension of such identities to the
metric $\rho$?

Close inspection of this question shows that the problem is
intimately related with the action of group (2.27) preserving
the metric $\rho$. One can check that the group $H$ is
transitive on each sphere $\goth S(r)$ only in  the two special cases
of  $s=1$ or $n=1$. For arbitrary $n$ and $s$ spheres
$\goth S(r)$ split into finitely many $H$-orbits. It turns out
that weight enumerators associated with such $H$-orbits satisfy
remarkable relations for mutually dual codes and distributions.
Such relations should be treated as a proper extension of the
MacWilliams identitiees to the metric $\rho$. 
 Unfortunately, within the limits of the present
paper we can not discuss these questions in detail. Here we consider
only the special case of $n=1$. The author hopes to return to
this matter elsewhere. 

\medskip{\bf Theorem 4.4.} {\rm(i)} {\it For every linear distribution
$D\s\BQ(q^s)$ its box and weight enumerators $\ph(D;z)$ and
$\Cal W(D;z)$ are related by 
$$
\Cal W(D;z)=z^s(1-z)\ph\left(D;\frac1z\right)+\#\{D\}z^{s+1}. 
\tag4.25
$$

{\rm(ii)} For linear mutually dual distributions $D$ and
$D^\perp\s\BQ(q^s)$ one has the identity 
$$
(qz-1)\Cal W(D;z)+1-z
$$
$$
=
\#\{D\}z^{s+1}\left(q(1-z)\Cal W
\left(D^\perp;\frac1{qz}\right)+qz-1\right),
\tag4.26
$$
or more symmetrically,
$$
v(D;z)=\#\{D\}q\,z^{s+2}v\left(D^\perp,\frac1{qz}\right),
\tag4.27
$$
where $v(D;z)=(qz-1)\Cal W(D;z)+1-z$. 
}
\medskip

{\it Proof.} By Lemma 3.1 every sphere $\goth S(r)\s\BQ(q^s)$
consists of a single fragment $g_r$. Precisely, the spheres have
the following structure: $\goth S(0)=[0,q^{-s})=\Dl^0_s$, and 
$$
\goth S(r)=[q^{-s+r-1},q^{-s+r})=\Dl^0_{s-r}\setminus
\Dl^0_{s-r+1}, 1\le r\le s, 
$$
where $\Dl^0_a=[0,q^{-a})$, $0\le a\le s$, are elementary boxes
(1.5) for $n=1$. Notice also that $\#\{D\cap\Dl^0_0\}=\#\{D\}$,
and $\#\{D\cap \Dl^0_s\}=1$, since $D$ is a linear distribution.

Substituting the foregoing relations to (4.24), we find that 
$$
\Cal
W(D,z)=1+\sum^s_{r=1}(\#\{D\cap\Dl^0_{s-r}\}-\#\{D\cap\Dl^0_{s-r+1}\})
z^r
$$
$$
=1+\sum^s_{a=0}\#\{D\cap\Dl^0_a\}z^{s-a}
-\sum^s_{a=0}\#\{D\cap\Dl^0_a\}z^{s+1-a}
$$
$$
-\#\{D\cap\Dl^0_s\}+\#\{D\cap\Dl^0_0\}z^{s+1}
$$
$$
=
\ph\left(D;\frac1z\right)z^s-\ph\left(D;\frac1z\right)z^{s+1}+
\#\{D\}z^{s+1}.
$$
This proves  relation (4.25). 

(ii) Identities (4.26) and (4.27) follow from  relations
(4.23) and (4.25) by a direct calculation. $\square$
\medskip

{\bf 4.5. Parity-check Matrices and $\rho$-Weights.} Linear
$[ns,k]_s$-codes can be treated as null subspaces of linear
surjections $\Mat_{n,s}(\BF_q)\to\BF^k_q$. For a given code $\CC$
such a mapping, written in a fixed basis, is said to be a
parity-check matrix of $\CC$. 

More precisely, let $\CH$ be an $ns$ by $k$ matrix over $\BF_q$
of rank $k\ge1$. Write $\CH=(\CH_1,\ldots,\CH_n)$, where
$\CH_j$, $1\le j\le n$, are $s$ by $k$ matrices. Then, a linear
$[ns,k]_s$-code $\CC \s\Mat_{n,s}(\BF_q)$ with the parity-check
matrix $\CH$ is defined by the following $k$ linear independent
equations 
$$
\CH_1\om^T+\ldots+\CH_n\om^T_n=0,
\tag4.28
$$
where $\Om=(\om_1,\ldots,\om_n)\in\CC$ and
$\om_j\in\Mat_{1,s}(\BF_q)$, $1\le j\le n$. 

We wish to evaluate the $\rho$-weight of a linear code in terms
of its parity-check matrix. For the Hamming weight the
corresponding result is well known (cf. [16, Chapter~1,
Theorem~10]). 

Introduce the following characteristic of a matrix
$\CH=(\CH_1,\ldots,\CH_n)$. For given integers $0\le d_j\le s$,
$1\le j\le n$, we form a set $Q(d_1,\ldots,d_n)$
of $d_1+\ldots+d_n$ vectors in $\BF^k_q$ selecting first $d_1$
columns from the matrix $\CH_1$, first $d_2$ columns from the
matrix $\CH_2$, etc. Now we put 
$$
\rho^\sharp(\CH)=\min\sum^n_{j=1} d_j,
\tag4.29
$$
where the minimum is extended over all integers
$(d_1,\ldots,d_n)\ne(0,\ldots,0)$ with linear {\it dependent}
sets $Q(d_1,\ldots,d_n)$. 

Characteristic (4.29) was originaly introduced in [20,
Definitions~6.8 and 7.1] for linear $(\dl,s,n)$-nets. It turns
out that quantity (4.29) coinsides with the $\rho$-weight of the
corresponding code (4.28). 

\medskip
{\bf Proposition 4.2.} {\it Let $\CC\s\Mat_{n,s}(\BF_q)$ be a
linear code with a parity-check matrix $\CH$. Then
$$
\rho(\CC)=\rho^\sharp(\CH).
\tag4.30
$$
}

\medskip
{\it Proof.} For each $\Om\in\CC\setminus\{0\}$ columns of the
matrices $\CH_j$, $1\le j\le n$, are related by (4.28).
 Comparing definitions (2.12), (2.14)  with (4.28),
we conclude the following: for given integers $0\le d_j\le s$,
$1\le j\le n$, $(d_1,\ldots, d_n)\ne(0,\ldots,0)$, there exists
a code word $\Om=(\om_1,\ldots,\om_n)\in\CC\setminus\{0\}$ with
$\rho(\om_j)=d_j$, $1\le j\le n$, if and only if the set
$Q(d_1,\ldots,d_n)$ is linear dependent. Since
$\rho(\Om)=d_1+\ldots+d_n$, equality (4.30) follows at once from
(4.29) and definition (4.1) for the weight $\rho(\CC)$.
$\square$ 
\medskip

\centerline{\bf 5. Constructions of Linear MDS Codes and Optimum
Distributions}
\medskip

{\bf 5.1 Hermite Interpolations over Finite Fields.}
Let $\BF_q[z]$, $q=p^e$, denote the ring of polynomials 
$$
f(z)=\sum^{t-1}_{i=0}f_i \,z^i
\tag5.1
$$
with coefficients $f_i\in\BF_q$, $t=\deg f+1$. For a given
polynomial (5.1) its $j^{th}$ hyperderivative $\cd^j
f\in\BF_q[z]$ is defined by
$$
\cd^j f(z)=\sum^{t-1}_{i=0}\left(\matrix i\\
j\endmatrix\right) f_i\,z^{i-j},
\tag5.2
$$
(cf. [15, Section 6.4]). From here on we  shall write $\left(\matrix i\\
j\endmatrix\right)$ for binomial coefficients modulo $p$,
moreover, as it usually is, $\left(\matrix i\\
j\endmatrix\right)=0$ for $j>i$. 

In particulary,  definition (5.2) implies
$$
\cd^j(z-\be)^i=\left(\matrix i\\
j\endmatrix\right)(z-\be)^{i-j}
\tag5.3
$$
for each $\be\in\BF_q$. Using (5.2) and (5.3), it is easy to check
the following expansion 
$$
f(z)=\sum^{t-1}_{j=0}\cd^j f(\be)(z-\be)^j,
\tag5.4
$$
where $\cd^j f(\be)$ denotes the value of $\cd^j f(z)$ at
$\be\in\BF_q$. 

It should be noted that for $j<p$  hyperderivative (5.2) and
the usual $j^{th}$ formal derivative $f^{(j)}(z)$ are related by
$$
\cd^j f(z)=\frac 1{j!}\,f^{(j)}(z).
\tag5.5
$$
Certainly, in the ring of polynomials over a field of infinite
characteristic  relation (5.5) holds for  all $j\in\BN_0$. 

We can formally add a point $\infty$ to $\BF_q$. In so doing,
all arithmetic operations (except $0\cdot\infty$) may be defined
in the obvious way on the set $\BF_q\cup\{\infty\}$ of $q+1$
elements. For a polynomial (5.1) of degree $t-1$ we put by
definition 
$$
\cd^j f(\infty)=
\cases f_{t-1-j}\quad&\text{if}\quad 0\le j\le t-1,\\
0\quad&\text{if}\quad j>t-1.
\endcases
\tag5.6
$$
Obviously, $\cd^j f(\infty)=\cd^j\check f(0)$, where $\check
f(z)=z^{t-1}f\left(\frac1 z\right)=f_{t-1}+f_{t-2}z+\ldots+f_0
z^{t-1}$ is a polynomial receprocal to $f(z)$. 

Let $\BM^t\s\BF_q[z]$ denote the set of all polynomials (5.1) of
degree less than $t$, $t\in\BN_0$, then $\BM^t$ is a vector
space of dimension $t$ over $\BF_q$. Consider the following {\it
Hermite interpolation problem}: Find a polynomial $f\in\BM^t$
which satisfies the equations 
$$
\cd^j f(\be_i)=a^{(j)}_i,\quad j=0,\ldots,t_i-1,\quad 
1\le i\le l,
\tag5.7
$$
where $l$ distinct elements
$\be_1,\ldots,\be_l\in\BF_q\cup\{\infty\}$ (interpolation nodes) are
fixed, and integers $t_i\in\BN$ and coefficients
$a^{(j)}_i\in\BF_q$ are given. Certainly, we assume here that
$q\ge l-1$. 

\medskip{\bf Proposition 5.1.} {\rm(i)} {\it The Hermite interpolation
problem $(5.7)$ has a unique solution $f\in\BM^t$, provided that
$t_1+\ldots+t_l=t$. 

{\rm(ii)} The homogeneous Hermite interpolation problem $(5.7)$
with all coefficients $a^{(j)}_i=0$ has a unique solution $f=0$
in the space $\BM^t$, provided that $t_1+\ldots+t_l\ge t$. 
}
\medskip

{\it Proof.} (i) At first, we suppose that all interpolation
nodes $\be_i\in\BF_q$, $1\le i\le l$. Consider the polynomials
$$
r_i(z)=\sum^{t_i-1}_{j=0}a^{(j)}_i(z-\be_i)^j,\quad 
1\le i\le l.
\tag5.8
$$
From  relations (5.3), (5.4), and (5.8) we conclude that 
equations (5.7) are equivalent to the following congruences in
the ring $\BF_q[z]$
$$
f(z)=r_i(z)\mod(z-\be_i)^{t_i},\quad 1\le i\le l.
\tag5.9
$$

In (5.9) the polynomials $(z-\be_i)^{t_i}$, $1\le i\le l$, are
pairwise coprime. Therefore,  congruences (5.9) have a unique
solution $f\in\BM^t$, $t=t_1+\ldots+t_l$ by the Chinese
remainder theorem in the ring $\BF_q[z]$ (see [16, Section 10.9,
Theorem~5]). 

Now, let one of interpolation nodes, say $\be_l$, coinsides with
$\infty$. Then, using (5.6), we put 
$$
f(z)=g(z)+r_l(z),
\tag5.10
$$
where
$$
r_l(z)=\sum^{t-1}_{j=t-t_l}a^{(j)}_l z^j,
$$
and $g(z)$ is a polynomial of degree less than $t-t_l$. 

Substituting (5.10) to first $l-1$ equations (5.7), we find that
the polynomial $g\in\BM^{t-t_l}$ is a solution of the following
Hermite interpolation problem with $l-1$ nodes
$\be_1,\ldots,\be_{l-1}\in\BF_q$:
$$
\cd^j g(\be_i)=b^{(j)}_i,\quad j=0,\ldots,t_i-1,\quad
1\le i\le l-1,
$$
where $b^{(j)}_i=a^{(j)}_i-\cd^j r(\be_i)$, and
$t_1+\ldots+t_{l-1}=t-t_l$. This problem has a unique solution
by the above arguments. 

(ii) This statement follows at once from (i) because we may
choose integers $t'_i\in\BN$, $t'_i\le t_i$, to satisfy the
equality $t'_1+\ldots+t'_l=t$. $\square$

\medskip
{\it Remark 5.1.} The Hermite interpolations may be
treated as a confluent Lagrange interpolation problem (cf. [5,
Chapter 2, Section 11]). This leads to another proof of
Proposition~5.1. Substituting  relations (5.2) for a
polynomial $f\in\BM^t$ to (5.7), we obtain a system of
$t_1+\ldots+t_l$ linear equations for $t$ unknown coefficients
$f_i$. If $t_1+\ldots+t_l=t$, then the system has a unique
solution, since the corresponding $t$ by $t$ matrix turns out  to be
a confluent Vandermonde matrix, and Proposition~1.1 
follows. However, the above arguments based on the Chinese
remainder theorem seem to be the most natural even in the case
of the Hermite interpolations over the field of complex numbers.

\medskip
{\bf 5.1. Explicit Constructions.}
Henceforward we assume that $q\ge n-1$. Fix $n$ distinct
elements $\be_1,\ldots,\be_n\in\BF_q\cup \{\infty\}$, and define
a linear mapping $\Ga_{n,s}:\BF_q[z]\to\Mat_{n,s}(\BF_q)$ by
setting 
$$
\Ga_{n,s}:\BF_q[z]\ni f\to \Ga_{n,s}f=
(\cd^{s-1-j}f(\be_i))\in\Mat_{n,s}(\BF_q),
\tag5.11
$$
where $0\le j\le s-1$ and $1\le i\le n$.

Thus, for each polynomial $f\in \BF_q[z]$ the matrix
$$
\Ga_{n,s}f=(\om^{(1)}_f,\ldots,
\om^{(n)}_f)^T\in\Mat_{n,s}(\BF_q)
$$
 consists of $n$ rows of the
form 
$$
\om^{(i)}_f=(\cd^{s-1}f(\be_i),\ldots,\cd
f(\be_i),f(\be_i))\in\Mat_{1,s}(\BF_q) ,\quad
1\le i\le n.
\tag5.12
$$

\medskip{\bf Lemma 5.1.} {\it For $t\le ns$ the image
$\Ga_{n,s}\BM^t\s 
\Mat_{n,s}(\BF_q)$ of the vector space $\BM^t\s\BF_q[z]$ under
mapping $(5.11)$ is a subspace of dimension $t$. 
}
\medskip

{\it Proof.} It suffices to prove that
$\Ga_{n,s}\BM^{ns}=\Mat_{n,s}(\BF_q)$, since $\BM^t\s\BM^{ns}$
for $t\le ns$ and the cardinalities of $\BM^{ns}$ and
$\Mat_{n,s}(\BF_q)$ coinside. Let $f\in\BM^{ns}$ and
$\Ga_{n,s}f=0$, then from (5.11) we conclude that $f$ is a
solution of a homogeneous Hermite interpolation problem (5.7)
with $l=n$, $t=ns$, $t_1=\ldots=t_n=s$. Therefore $f=0$ by
Proposition~5.1. $\square$ 
\medskip

\medskip{\bf Theorem 5.1.} {\it For each integer $1\le k\le ns$ the
subspace $\Ga_{n,s}\BM^{k}\s\Mat_{n,s}(\BF_q)$ is an MDS
$[ns,k]_s$-code in the metric $\rho$. 
}
\medskip

{\it Proof.} By Definition 2.1 and  relations (4.1) it
sufficies to show that $\rho(\Ga_{n,s}f)\ge  ns-k+1$ for all
$f\in\BM^{ks}\setminus\{0\}$. Suppose on the contrary that there
is a non-zeroth polynomial $f\in\BM^{ks}$ such that
$\rho(\Ga_{n,s}f)\le  ns-k$. By  definitions (2.14), and
(5.11), (5.12) this yields 
$$
\rho(\Ga_{n,s}f)=\rho(\om^{(1)}_f)+\ldots+\rho(\om^{(n)}_f)\le ns-k.
\tag5.13
$$

From (5.13) we conclude that the strict inequality 
$$
\rho(\om^{(i)}_f)<s
\tag5.14
$$
holds at least for one of indices $i$. Without loss of
generality, we may assume that  inequality (5.14) holds for
$i=1,\ldots,l$ with $l\ge 1$ and $\rho(\om^{(i)}_f)=s$ for
$i=l+1,\ldots,n$. 

Using (5.11), (5.12), and (5.14), we find that the polynomial
$f\in\BM^{k}$ is a solution of the following homogeneous
Hermite interpolation problem (5.7) 
$$
\cd^j f(\be_i)=0,\quad j=0,\ldots,t_i-1,\quad 1\le i\le l,
$$
where $t_i=s-\rho(\om^{(i)}_f)\in\BN$. 

It is easy to check that $t_1+\ldots+t_l=ns-\rho(\Ga_{n,s}f)$.
Therefore, $t_1+\ldots+ t_l\ge k$ and $f=0$ by Proposition~1.1,
a contradiction. $\square$
\medskip

Theorem 5.1 gives broad classes of linear MDS $[ns,k]_s$-codes in
the metric $\rho$ with arbitrary values of the parameters $n,k$,
and $s$. By Theorem~2.1 such MDS codes correspond to linear
optimum $[ns,k]_s$-distributions, which can be explicitly
described as follows. Introduce a linear mapping
$\ga_{n,s}:\BF_q[z]\to \BQ^n(q^s)$ by setting 
$$
\ga_{n,s}:\BF_q[z]\ni f\to\ga_{n,s}f=
X=(x_1,\ldots,x_n)^T\in\BQ^n(q^s),
\tag5.15
$$
where coordinates $x_i$ of the point $X$ are given by 
$$
x_i=\sum^s_{j=1}\cd^{j-1} f(\be_i)q^{-j}.
\tag5.16
$$

Using (5.1), (5.2), and (5.6), we can write the coefficients
$\cd^{j-1}f(\be_i)$ in the form 
$$
\cd^{j-1}f(\be_i)=\sum^{t-1}_{l=0}\left(\matrix j-1\\
l\endmatrix\right) \be_i^{l+1-j}f_l
\tag5.17
$$
for $\be_i\in\BF_q$, and 
$$
\cd^{j-1}f(\be_i)=f_{t-j}
\tag5.18
$$
for $\be_i=\infty$.

\medskip{\bf Theorem 5.2.} {\it For each integer $1\le k\le ns$
the 
subspace $\ga_{n,s}\BM^{k}\s\BQ^n(q^s)$ is an optimum
$[ns,k]_s$-distribution. 
}
\medskip

{\it Proof.} It suffices to compare  relations (5.15), (5.16)
with (5.11), (5.12), and use Theorems~5.1 and 2.1. $\square$
\medskip

In the special case of $k=s$ an optimum $[ns,s]_s$-distribution
$\ga_{n,s}\BM^s$ (defined by the relations from (5.15) to
(5.18)) coinsides with a $(0,s,n)$-net of zeroth deficiency (cf.
Proposition~1.1). Nets of such a kind were
 given previously in [11] and [20]. 

\medskip
\centerline{\bf 6. Reconstructions of Codes and Distributions} 
\medskip

{\bf 6.1. Peano's Bijection.}
As usual, for a matrix $\Om\in\Mat_{g,s}(\BF_q)$ we write
$\Om=(\om_1,\ldots,\om_g)^T$ where
$\om_j=(\xi^{(j)}_1,\ldots,\xi^{(j)}_s)$ $\in \Mat_{1,s}(\BF_q)$.
Introduce a mapping $\pi_g$ from $\Mat_{g,s}(\BF_q)$ to
$\Mat_{1,gs}(\BF_q)$ by 
$$
\gathered
\pi_g:\Mat_{g,s}(\BF_q)\ni\Om\to\pi_g\Om
=
\\
=(\xi^{(1)}_{1},\ldots,\xi^{(g)}_1,\ldots,\xi^{(1)}_s,\ldots,\xi^{(g)}_s)
\in\Mat_{1,gs}(\BF_q).
\endgathered
\tag6.1
$$
It is obvious that  mapping (6.1) is an isomorphism of vector
spaces $\Mat_{g,s}(\BF_q)$ and  $\Mat_{1,gs}(\BF_q)$. 

By  relations (2.4), (2.5), and (2.6)  mapping (6.1) can
be considered as an isomorphism of vector spaces $\BQ^g(q^s)$
and $\BQ(q^{gs})$. In this case $\pi_g$ coinsides with the well
known Peano's mapping (restricted on $\BQ^g(q^s))$ which gives a
bijection between points of the unit cube $U^g$ and the unit
segment $[0,1)$. 

More generally, we introduce an isomorphism $\pi_{g,n}$ of
vector spaces $\Mat_{gn,s}(\BF_q)$ and  $\Mat_{n,gs}(\BF_q)$ as
follows. For a matrix $\Om\!\!\in\!\!\Mat_{gn,s}(\BF_q)$ we write
$\Om\!\!=\!\!(\Om_1,...,\Om_n)^T$ where $\Om_j\in\Mat_{g,s}(\BF_q)$,
$1\le j\le n$. Now we put 
$$
\pi_{g,n}:\Mat_{gn,s}(\BF_q)\ni\Om\to \pi_{g,n}\Om
$$
$$
=(\pi_g\Om_1,\ldots,\pi_g\Om_n)^T\in\Mat_{n,gs}(\BF_q).
\tag6.2
$$

In other words, we split $gn$ rows of a matrix
$\Om\in\Mat_{gn,s}(\BF_q)$ into $n$ subsets
$\Om_1,\ldots,\Om_n$, each of which consists of successive $g$
rows; and subsequently,  we map by $\pi_g$ first $g$ rows of $\Om_1$
to the first row of the matrix
$\pi_{g,n}\Om\in\Mat_{n,gs}(\BF_q)$, second $g$ rows of $\Om_2$ to
the second row of $\pi_{g,n}\Om$, and so on, last $g$ rows
$\Om_n$ are mapped to the $n^{th}$ row of $\pi_{g,n}\Om$. 

Using  relations (2.8) and (6.2), we define an isomorphism 
$$
\pi_{g,n}:\BQ^{gn}(q^s)\ni X\to\pi_{g,n}X\in \BQ^n(q^{gs})
\tag6.3
$$
of vector spaces $\BQ^{gn}(q^s)$ and $\BQ^n(q^{gs})$  by setting
$$
\Om\l\pi_{g,n}X\r=\pi_{g,n}\Om\l X\r.
\tag6.4
$$

It is obvious that  mapping (6.3), (6.4) is a restriction on
$\BQ^{gn}(q^s)$ and $\BQ^n(q^{gs})$ of the well known Peano's
bijection between points of the unit cubes $U^{gn}$ and $U^n$. 

The foregoing implies the following.

\medskip{\bf Lemma 6.1.} {\rm(i)} {\it Let
$\CC\s\Mat_{gn,s}(\BF_q)$ be 
a $[gns,gk]_s$-code and $D\s\BQ^{gn}(q^s)$ be a
$[gns,gk]_s$-distribution, then
$\pi_{g,n}\CC\s\Mat_{n,gs}(\BF_q)$ is a $[ngs,gk]_{gs}$-code and
$\pi_{g,n}D\s\BQ^n(q^{gs})$ is a $[ngs,gk]_{gs}$-distribution. 

{\rm(ii)} If a code $\CC\s\Mat_{gn,s}(\BF_q)$ and a distribution
$D\s\BQ^{gn}(q^s)$ are linear, then the code
$\pi_{g,n}\s\Mat_{n,gs}(\BF_q)$ and the distribution $\pi_{g,n}
D\s\BQ^n(q^{gs})$ are also linear.
}
\medskip

Now our concern is with the behavior of the Hamming and
non-Hamming weights $\ka$ and $\rho$ under  mappings (6.2)
and (6.3). 

\medskip{\bf Lemma 6.2.} {\rm(i)} {\it For arbitrary matrices
$\Om\in\Mat_{gn,s}(\BF_q)$ and points $X\in\BQ^{gn}(q^s)$ one
has the relations 
$$
\ka(\pi_{g,n}\Om)=\ka(\Om),\quad
\ka(\pi_{g,n}X)=\ka(X),
\tag6.5
$$
and
$$
\rho(\pi_{g,n}\Om)\ge\rho(\Om),\quad
\rho(\pi_{g,n}X)\ge\rho(X).
\tag6.6
$$

{\rm(ii)} For arbitrary codes $\CC\s\Mat_{gn,s}(\BF_q)$ and
distributions $D\s\BQ^{gn}(q^s)$ one has the relations
$$
\ka(\pi_{g,n}\CC)=\ka(\CC),\quad
\ka(\pi_{g,n}D)=\ka(D)
\tag6.7
$$
and
$$
\rho(\pi_{g,n}\CC)\ge\rho(\CC),\quad
\rho(\pi_{g,n}D)\ge\rho(D).
\tag6.8
$$
}
\medskip

{\it Proof.} (i) It suffices to prove (6.5) and (6.6) only for
matrices $\Om\in\Mat_{gn,s}(\BF_q)$. 

To prove (6.5) we simply notice that the matrix $\pi_{g,n}\Om$
is obtained by a reordering (see (6.1), (6.2)) of elements of
$\Om$. Therefore, the matrices $\pi_{g,n}\Om$ and $\Om$ have the
same number of non-zero entries.

To prove (6.6) we first consider the behavior of the
$\rho$-weight under  mapping (6.1). We wish to prove that for
an arbitrary matrix $\Om\in\Mat_{g,s}(\BF_q)$ one has the
inequality 
$$
\rho(\pi_g\Om)\ge\rho(\Om).
\tag6.9
$$

Since (6.9) holds trivially for $\Om=0$, we suppose that
$\Om\ne0$. Let $\Om=(\om_1,\ldots,\om_g)^T$ where
$\om_j\in\Mat_{1,s}(\BF_q)$, $1\le j\le g$. Introduce an integer
$l$, $1\le l\le g$, by setting 
$$
l=l(\Om)=\max\{j:\rho(\om_j)\ne0\}.
\tag6.10
$$

Comparing  definitions (2.12) and (6.1), we conclude that 
$$
\rho(\pi_g\Om)=\rho(\om_l)+(l-1)s.
\tag6.11
$$
Moreover,  definition (2.14) implies 
$$
\rho(\Om)=\sum^l_{j=1}\rho(\om_j).
\tag6.12
$$

Using  relations (2.13), (6.11), and (6.12), we find that 
$$
\rho(\Om)=\rho(\om_l)+\sum^{l-1}_{j=1}\rho(\om_j)\le
\rho(\om_l)+(l-1)s=\rho(\pi_g\Om).
$$

The proof of  inequality (6.9) is complete. 

Now, let an arbitrary matrix $\Om\!\in\!\Mat_{gn,s}(\BF_q)$ be
given. Then $\Om\!=\!(\Om_1,...,\Om_n)^T$, where
$\Om_j\in\Mat_{g,s}(\BF_q)$, $1\le j\le n$. Comparing 
definitions (2.14) and (6.2), we conclude that 
$$
\rho(\Om)=\sum^n_{j=1}\rho(\Om_j),\quad
\rho(\pi_{g,n}\Om)=\sum^n_{j=1}\rho(\pi_g\Om_j).
\tag6.13
$$

Substituting  inequality (6.9) to (6.13), we obtain 
inequality (6.6). 

(ii) The proposition (ii) follows at once from (i) and 
definitions (2.18), (2.19) of weights $\ka$ and $\rho$.
$\square$ 
\medskip

\medskip{\bf Proposition 6.1.} {\it Let
$\CC\s\Mat_{gn,s}(\BF_q)$ be an 
MDS $[gns,gk]_s$-code in the metric $\rho$ and $D\s\BQ^{gn}(q^s)$
be an optimum $[gns,gk]_s$-distribution. Suppose that $k=st$,
$1\le t\le n-1$. Then, 
$$
\rho(\pi_{g,n}\CC)=(n-t)gs+1,\quad
\rho(\pi_{g,n}D)=(n-t)gs+1.
\tag6.14
$$
Hence, $\pi_{g,n}\CC\s\Mat_{n,gs}(\BF_q)$ is an MDS
$[ngs,gk]_{gs}$-code in the metric $\rho$ and
$\pi_{g,n}D\s\BQ^n(g^{gs})$ is an optimum
$[ngs,gk]_{gs}$-distribution. 

Furthermore, the corresponding Hamming weights satisfy the
following inequalities 
$$
\ka(\pi_{g,n}\CC)\ge(n-t)g+1,\quad
\ka(\pi_{g,n}D)\ge (n-t)g+1.
\tag6.15
$$
}
\medskip

{\it Proof.} It suffices to prove  relations (6.14) and
(6.15) only for codes. 

To prove (6.14) we notice that $\rho(\pi_{g,n}\CC)\ge(n-t)gs+1$
by Definition~2.1 and  inequality (6.8). It remains to refer
to Proposition~2.1.

To prove (6.15) we notice that 
$$
\ka(\CC)\ge s^{-1}\rho(\CC)=(n-t)g+s^{-1}
$$
by  inequality (2.15) and Definition~2.1. Therefore
$\ka(\CC)\ge(n-t)g+1$, since the weight $\ka(\CC)$ is an
integer. It remains to refer to  equality (6.7). $\square$ 
\medskip

{\bf 6.2. Duality and Peano's Bijection.}
We wish to consider the behavior of  inner products (4.2) and
(4.3) under  mappings (6.1) and (6.2). We introduce an
isomorphism $J_g$ of the space $\Mat_{g,s}(\BF_q)$ by setting 
$$
\gathered
J_g:\Mat_{g,s}(\BF_q)\ni \Om=(\om_1,\ldots,\om_g)^T\to J_g\Om
\\
=
(\om_g,\ldots,\om_1)^T\in\Mat_{g,s}(\BF_q),
\endgathered
\tag6.16
$$
where $\om_j\in\Mat_{1,s}(\BF_q)$, $1\le j\le g$. Thus, $J_g$ is
merely a permutation of rows in inverse order. 

By a direct calculation from  definitions (4.2), (4.3) and
(6.1), (6.16) we derive the identity 
$$
\l\pi_g\Om_1,\pi_g\Om_2\r=\l J_g\Om_1,\Om_2\r=
\l\Om_1,J_g\Om_2\r
\tag6.17
$$
for all $\Om_1$, $\Om_2\in\Mat_{g,s}(\BF_q)$. 

Now we introduce an isomorphism $J_{g,n}$ of the space
$\Mat_{gn,s}(\BF_q)$ by setting
$$
J_{g,n}:\Mat_{gn,s}(\BF_q)\ni \Om=
(\Om_1,\ldots,\Om_n)^T\to J_{g,n}\Om
$$
$$
=(J_g\Om_1,\ldots,J_g\Om_n)^T\in \Mat_{gn,s}(\Bbb F_q),
\tag6.18
$$
where $\Om_j\in\Mat_{g,s}(\BF_q)$.

Using  definitions (4.2), (4.3), (6.1), (6.2), and 
relation (6.17), we obtain the identity
$$
\l\pi_{g,n}\Om_1,\pi_{g,n}\Om_2\r=\l J_{g,n}\Om_1,\Om_2\r
=\l\Om_1,J_{g,n}\Om_2\r
\tag6.19
$$
for all $\Om_1$, $\Om_2\in\Mat_{gn,s}(\BF_q)$.

\medskip{\bf Proposition 6.2.} {\it Given a linear code
$\CC\s\Mat_{gn,s}(\BF_q)$ and its image
$\pi_{g,n}\CC\s\Mat_{n,gs}(\BF_q)$ under  mapping $(6.2)$.
Then the dual codes $\CC^\perp$ and $(\pi_{g,n}\CC)^\perp$ are
related by 
$$
(\pi_{g,n}\CC)^\perp=\pi_{g,n}(J_{g,n}\CC^\perp).
\tag6.20
$$
}
\medskip

{\it Proof.} Since  mapping (6.2) is a bijection, we can
write $(\pi_{g,n}\CC)^\perp=\pi_{g,n}\CC_1$ with a linear code
$\CC_1\s\Mat_{n,gs}(\BF_q)$. Using Definition~4.1 and  relation
(6.19), we conclude that $J_{g,n}\CC_1=\CC^\perp$. Therefore,
$\CC_1=J_{g,n}\CC^\perp$, because $J^2_{g,n}$ is the identity,
and  relation (6.20) follows. $\square$
\medskip

For linear codes and distributions Proposition 6.1 can be
supplemented with the following result. 

\medskip{\bf Proposition 6.3.} {\it Let
$\CC\s\Mat_{gn,s}(\BF_q)$ be a 
linear MDS $[gns,gk]_s$-code in the metric $\rho$ and
$D\s\BQ^{gn}(q^s)$ be a linear optimum
$[gns,gk]_s$-distribution. Suppose that $k=st$, $1\le t\le n-1$.
Then, the Hamming and non-Hamming weights of the dual subjects
$(\pi_{g,n}\CC)^\perp$ and $(\pi_{g,n}D)^\perp$ satisfy the
following relations 
$$
\rho((\pi_{g,n}\CC)^\perp)=kg+1,\quad
\rho((\pi_{g,n}D)^\perp)=kg+1
$$
and 
$$
\ka((\pi_{g,n}\CC)^\perp)\ge tg+1,\quad
\ka((\pi_{g,n}D)^\perp)\ge tg+1.
$$
}
\medskip

{\it Proof.} It suffices to prove Proposition 6.3 only for
codes. First of all, we notice that 
$$
\rho(J_{g,n}\CC_1)=\rho(\CC_1)
\quad
\ka(J_{g,n}\CC_1)=\ka(\CC_1)
\tag6.21
$$
for an arbitrary code $\CC\s\Mat_{gn,s}(\BF_q)$, because the
mapping $J_{g,n}$ is nothing but a permutation of rows of the
corresponding code words, that preserves both weights $\rho$ and
$\ka$ (cf. Remark~2.2). 

Now let $\CC\s\Mat_{gn,s}(\BF_q)$ be a linear MDS
$[gns,gk]_s$-code in the metric $\rho$. Then, by Theorem~4.1 (i),
$\CC^\perp\s \Mat_{gn,s}(\BF_q)$ is an MDS $[gn,g(n-k),s]$-code in
the metric $\rho$. Therefore,  relations (6.7), (6.15) (where
$k$ is replaced by $ns-k$), (6.20), and (6.21) imply 
$$
\rho((\pi_{g,n}C)^\perp)=\rho(\pi_{g,n}(J_{g,n}C^\perp))=
\rho(J_{g,n}C^\perp)=\rho(C^\perp)=kg+1,
$$
and similarly, 
$$
\ka((\pi_{g,n}C)^\perp)=\ka(\pi_{g,n}(J_{g,n}C^\perp))=
\ka(J_{g,n}C^\perp)=\ka(C^\perp)\ge tg+1. \ \ \square
$$
\medskip

{\bf 6.3. Codes and Distributions with large $\rho$- and
$\ka$-Weights.}
Now we can give explicit constructions of codes and
distributions with large weights simultaneously in both metrics
$\rho$ and $\ka$. 

Given arbitrary integers $g\!\in\!\BN$ and $k=st$, $1\le t\le
n-1$.  Define a linear
$[gns,gk]_s$-code $\CC^{(g)}\!\s\!\Mat_{gn,s}(\BF_q)$ and a linear
$[gns,gk]_s$-distribution $D^{(g)}\s\BQ^{gn}(q^s)$ by 
$$
\CC^{(g)}=\Ga_{gn,s}\BM^{gk},\quad
D^{(g)}=\ga_{gn,s}\BM^{gk},
\tag6.22
$$
where the mappings $\Ga_{gn,s}$ and $\ga_{gn,s}$ are given in
(5.11). Certainly, we assume here that 
$$
q=p^e\ge gn-1.
\tag6.23
$$

Next, we introduce an $[ngs,gk]_{gs}$-code
$\pi_{g,n}\CC^{(g)}\s\Mat_{n,gs}(\BF_q)$ and an
$[ngs,gk]_{gs}$ dis\-tri\-bu\-ti\-on
$\pi_{g,n}D^{(g)}\s\BQ^n(q^{gs})$, where 
$\pi_{g,n}$ are Peano's bijections (6.2) and (6.3). We consider 
also the dual $[ngs,ngs-k]_{gs}$-code
$(\pi_{g,n}\CC^{(g)})^\perp\s\Mat_{n,gs}(\BF_q)$ and
$[ngs,ngs-k]_{gs}$-distribution
$(\pi_{g,n}D^{(g)})^\perp\s\BQ^n(q^{gs})$.  

It is worthwhile writing coordinates of points
$X=(x_1,\ldots,x_n)\in D^{(g)}$ explicitly. Comparing relations
(5.15), (5.16) with (6.1), (6.2) and (6.22), we find that 
$$
x_j=\sum^g_{l=1}q^{-(l-1)s}\sum^s_{i=1} \cd^{i-1}f(\be_{j,l})
q^{-i}, 
\tag6.24
$$
where $f\in\Bbb M^{gk}$, and $\be_{j,l}$, $1\le j\le n$, $1\le
l\le g$, are $ng$ distinct elements in $\BF_q\cup\{\infty\}$.
Notice that such elements $\be_{j,l}$ exist by (6.23).

\medskip{\bf Theorem 6.1.} {\it  With the above notation one has the
following relations for the corresponding Hamming and
non-Hamming weights 
$$
\gathered
\rho(\pi_{g,n}\CC^{(g)})=(ns-k)g+1,\quad 
\rho(\pi_{g,n}D^{(g)})=(ns-k)g+1,\\
\ka(\pi_{g,n}\CC^{(g)})\ge(n-t)g+1,\quad
\ka(\pi_{g,n}D^{(g)})\ge(n-t)g+1,
\endgathered
\tag6.25
$$
and
$$
\gathered
\rho((\pi_{g,n}\CC^{(g)})^\perp)=kg+1,\quad 
\rho((\pi_{g,n}D^{(g)})^\perp)=kg+1,\\
\ka((\pi_{g,n}\CC^{(g)})^\perp)\ge tg+1,\quad
\ka((\pi_{g,n}D^{(g)})^\perp)\ge tg+1,
\endgathered
\tag6.26
$$
}
\medskip

{\it Proof.} As usual, we prove Theorem~6.1 only for codes.  
By Theorem 5.1 the code $C^{(g)}\s\Mat_{gn,s}(\BF_q)$ given in
(6.22)  is a
linear MDS $[gns,gk]_s$-code in the metric $\rho$. This implies
 relations (6.25) by Proposition 6.1 and  relations (6.26)
by Proposition~6.3. $\square$
\medskip

{\bf 6.4. Variations of the Base $q$.}
Every distribution in the space $\BQ^n(q^s)$, $q=p^e$, can be
treated as a distribution in the space $\BQ^n(p^{es})$ (cf.
Lemma~1.1). We wish to find out what happens to metric
properties of a given distribution $D\s\BQ^n(q^s)$ when the base
$q$ varies from $p^e$ to $p$. 

In what follows we write $\ka_q$, $\rho_q$ for the corresponding
Hamming and non-Hamming weights on the space $\BQ^n(q^s)$,
$q=p^e$, and $\ka_p$, $\rho_p$ for those on the space
$\BQ^n(p^{es})$ (cf. Sec.~2.2). 

Let $x\in\BQ(q^s)$. Consider  representation (2.5) in bases
$q$ and $p$, 
$$
x=\sum^s_{i=1}\xi_i(x)q^{i-s-1}=
\sum^{es}_{j=1}\th_j(x)p^{j-es-1},
\tag6.27
$$
where
$\xi_i(x)=\xi_i=\mu^{(1)}_i+\mu^{(2)}_ip+\ldots+\mu^{(e)}_ip^{e-1}$,
and $\mu^{(e)}_i=\mu^{(e)}_i(x)\in\{0,1,\ldots,p-1\}$ (cf.
(2.6)). The coefficients $\th_j=\th_j(x)$ in (6.27) are
determined by the following rule for the corresponding code
words: 
$$
\Mat_{1,s}(\BF_q)\ni \om_q\l x\r=(\xi_1,\ldots,\xi_s)
$$
$$
=
(\mu^{(1)}_1,\ldots,\mu^{(e)}_1,\mu^{(1)}_2,\ldots,\mu^{(e)}_2,\ldots,
\mu^{(1)}_s,\ldots,\mu^{(e)}_s)
$$
$$
=(\th_1,\ldots,\th_{es})=\om_p\l x\r\in\Mat_{1,es}(\BF_p).
\tag6.28
$$

Since $\rho_q(0)=\rho_p(0)=0$, we may  consider only 
$x>0$. With 
the above notation  definition (2.12) takes the form 
$$
\rho_q(x)=\max\{i:\xi_i\ne0\},\quad
\rho_p(x)=\max\{j:\th_j\ne0\}.
\tag6.29
$$

Suppose that $\rho_q(x)=r$, $r\in\BN$, then using (6.28) and
(6.29), we conclude that $\th_j=0$ for all $j>er$ and there is
an element $\th_j\ne0$ with $e(r-1)+1\le j\le er$. Therefore, 
$$
e(\rho_q(x)-1)+1\le\rho_p(x)\le e\rho_q(x).
\tag6.30
$$

Combining inequality (6.30) with  definition (2.14), we
find the following relation
$$
e(\rho_q(X)-1)+1-(e-1)(n-1)\le\rho_p(X)\le e\rho_q(X)
\tag6.31
$$
where $X\in\BQ^n(q^s)$ is an arbitrary point. Certainly, the
left side of (6.30) and (6.31) is informative only for
sufficiently large weights $\rho_q$. 

Similarly, using  rule (6.28), we obtain the following
inequality for the Hamming weights 
$$
\ka_q(X)\le\ka_p(X)\le e\ka_q(X),
\tag6.32
$$
where $X\in\BQ^n(q^s)$ is an arbitrary point. Notice that
bounds of such a kind are well known in coding theory (cf.
[16, Section 7.7]). 

The foregoing immediately implies the following result.

\medskip{\bf Proposition 6.4.} {\it Given a distribution $D\s\BQ^n(q^s)$,
$q=p^e$. Then, its Hamming and non-Hamming weights 
in bases $q$ and $p$ are related by 
$$
e(\rho_q(D)-1)+1-(e-1)(n-1)\le \rho_p(D)\le e\rho_q(D)
\tag6.33
$$
and 
$$
\ka_q(D)\le\ka_p(D)\le e\ka_q(D).
\tag6.34
$$

In particular, if $D\s\BQ^n(q^s)$ is an optimum
$[ns,k]_s$-distribution in base $q$ then its weight $\rho_p(D)$
in base $p$ satisfies the inequality 
$$
\rho_p(D)\ge (ns-k)e+1-(e-1)(n-1). 
\tag6.35
$$
}
\medskip

{\it Proof.} Inequalities (6.33), (6.34) follow from (6.31),
(6.32), and  definition (2.18). If $D\s\BQ^n(q^s)$ is an
optimum $[ns,k]_s$-distribution, then $\rho_q(D)= ns-k+1$ by
Theorem~2.1, and  inequality (6.35) follows from (6.33).
$\square$ 
\medskip

We return to the $[ns,k]_s$-distribution $\pi_{g,n}D^{(g)}\s
\BQ^n(q^{gs})$  and the dual $[ns,ns-k]_s$-distribution
$(\pi_{g,n}D^{(g)})^\perp\s\BQ^n (q^{gs})$ which were introduced in
Sec.~6.3. The bounds for their Hamming and non-Hamming weights
$\ka=\ka_q$ and $\rho=\rho_q$ in base $q$ are given in
Theorem~6.1. Now we wish to estimate the corresponding weights
$\ka_p$ and $\rho_p$ in base $p$. 

\medskip{\bf Theorem 6.2.} {\it With the notation of Sec.~$6.3$
one has 
the following relations for the corresponding Hamming and
non-Hamming weights in base $p$
$$
\gathered
\rho_p(\pi_{g,n}D^{(g)})\ge (ns-k)eg+1-(e-1)(n-1),\\
\ka_p(\pi_{g,n}D^{(g)})\ge(n-t)g+1,
\endgathered
\tag6.36
$$
and
$$
\gathered
\rho_p((\pi_{g,n}D^{(g)})^\perp)\ge k\,eg+1-(e-1)(n-1),\\
\ka_p((\pi_{g,n}D^{(g)})^\perp)\ge t \,eg+1,
\endgathered
\tag6.37
$$
where integers $e$, $g\in\BN$ and the prime $p$ are related by
$(6.23)$. 
}
\medskip

{\it Proof.} It suffices to compare Theorem~6.1 and
Proposition~6.4; in so doing  inequalities (6.36) and (6.37)
follow from  inequalities (6.25) and (6.26), respectively.
$\square$
\medskip

Thus, we have constructed explicitly the distributions
$\pi_{g,n}D^{(g)}$ and $(\pi_{g,n}D^{(g)})^\perp$ with large weights
simulteneously in metrics $\rho_p$ and $\ka_p$. As mentioned in
the Introduction distributions with such metric properties are of
crucial importance for further applications (see [7]). 

{\advance \parindent by 2mm
\Refs
\item{[1]}  Adams MJ, Shader~BL
(1997)  A construction for $(t,m,s)$-nets in base $q$.
SIAM J Discrete Math {\bf 10}: 460--468
\item{[2]}   Baker RC 
  On irregularities of
distribution. II.  J London Math Soc
(to appear)
\item{[3]} Beck J (1989) A two-dimensional van
Aardenne--Ehrenfest theorem in irregularities of distributions.
Composito Math {\bf 72}: 269--339
\item{[4]} Beck J, Chen WWL (1987)
 Irregularities of Distributions. Cambridge: Univ Press,
Cambridge 
\item{[5]} Berezin IS, Zhidkov NP (1965) Computing Methodes.
{\bf1}: Addison--Wesley, Reading Mass
\item{[6]} Chen WWL (1983) On irregularities of distributions. II.
 Quart J Math Oxford (2). {\bf34}: 257--279
\item{[7]} Chen WWL, Skriganov MM (1999)
Explicit constructions of point sets with
a minimal order of the $L_2$-discrepancy. Macquarie Univ \&
Steklov Math Inst St. Petersburg Depart. Preprint 
\item{[8]} Clayman AT, Mullen GL (1997) 
Improved $(t,m,s)$-net parameters from the
Gilbert--Var\-sha\-mov bound.
Applicable Algebra Engrg Comm Comp {\bf 8}: 491--496
\item{[9]} Conway JH, Sloane NJ (1988)  Sphere
Packings, Lattices and Groups.  New York: Springer
\item{[10]} Edel Y, Bierbrauer J (1997)
Construction of digital nets from BCH-codes.
In: H.~Niederreiter et al. (eds)
 Monte Carlo and Quasi -- Monte Carlo Methods
1996. Lecture
Notes in Statistics {\bf127}, pp.~221--231
New York:  Springer
\item{[11]} Faure H (1982)
Descr\'epance de suites associ\'ees 
\`a un syst\`me de num\'eration (en dimension $s$).
 Acta Arith  {\bf 41}: 337--351
\item{[12]} Kuipers L, Niederreiter H (1974)
Uniform Distribution of
Sequences.   New York: John Wiley \& Sons etc
\item{[13]} Lawrence KM, Mahalanabis A, Mullen GL,
Schmid WCh (1996) Construction of digital $(t,m,s)$-nets
from linear codes: In: Cohen S and Niederreiter H (eds)
  {Finite Fields and Applications
 London Math Soc
Lecture Note Series} {\bf233}: pp.~189--208,
Cambridge: Cambridge Univ Press
\item{[14]} Lang S (1971) Algebra, Addison--Wesley,
Reading, Mass
\item{[15]} Lidl R, Niederreiter H (1983) Finite Fields.
Addison--Wesley, London etc
\item{[16]} MacWilliams FJ, Sloane NJA (1977)
The Theory of error-correcting codes, North-Holland,
Amsterdam
\item{[17]} Mandelbaum DM (1979)
Construction of error-correcting
codes by interpolation.   IEEE Trans Inf Theor
{\bf IT-25}: 27--35
\item{[18]} Martin WJ, Stinson DR Association schemes for
ordered orthogonal arrays and $(T,M,S)$-nets. Cannad J Math (to
appear) 
\item{[19]} Matou\v sek J (1998)
 Geometric Discrepancy (an Illustrated
Guide).  Berlin: Springer
\item{[20]} Niederreiter H (1987) Point sets and sequences with 
small discrepancy.  Monatsh Math {\bf 104}: 273--337
\item{[21]} Niederreiter H (1992)  Random Number
Generation and Quasi -- Monte Carlo Methods.
SIAM, Philadelphia
\item{[22]} Niederreiter H (1998) Nets, $(t,s)$-sequences, and
algebraic curves over finite fields with many rational
points.   International Congress of Mathematicians
in Berlin, Documenta Mathematica, Extra Volume III:
 377--386
\item{[23]} Niederreiter H, Xing CP (1996)
Low-discrepancy sequences and global function fields with 
many rational places.  Finite Fields Appl
{\bf 2}:  241--273
\item{[24]} Niederreiter H, Xing CP (1998) Nets, $(t,s)$-sequences,
and algebraic geometry. In: Hellekalek P and Larcher G (eds) Pseudo- and
Quasi-Random Point Sets 
Lectures Notes in Statistics, {\bf 138},  New York:  Springer
\item{[25]} Rosenbloom MYu, Tsfasman MA (1997) Codes in the
$m$-metric. Problem i Peredachi Inf {\bf33} No~1: 55--63 (in
Russian); English transl in: Problems of Inf Trans {\bf33}
No~1: 45--52 (1997) 
\item{[26]} Skriganov MM (1994) Constructions of uniform
distributions in terms of geometry of numbers.
Algebra i Analiz {\bf 6}, No~3: 200--230;
 Reprinted in: St.~Petersburg Math J, {\bf6}:
635--664 (1995)
\item{[27]} Skriganov MM (1998) Ergodic theory on $SL(n)$,
diophantine approximations, and anomalies in the lattice
point problem.  Invent Math {\bf132}:  1--72
\item{[28]} Skriganov MM (1998) Uniform distributions, error-correcting
codes, and interpolations over finite fields.
Steklov Math Inst St. Petersburg Depart. Preprint 
\item{[29]} Sobol IM (1967) On the distribution of points in a cube 
and the approximate evaluation of integrals. \v Z~Vy\v cisl
Mat i Mat Fiz {\bf7}: 784--802; (in Russian) English transl in
USSR Comp Math and Math Phys {\bf 7}:
86--112
\item{[30]} Sobol IM (1969)  Multi-Dimensional Quadrature
Formulas and Haar Functions. Nauka, Mos\-cow 
(in Russian)
\item{[31]} Stichtenoth H (1993) Algebraic Function Fields 
and Codes.   Berlin: Springer
\item{[32]} Tsfasman MA, Vl\v adut SG (1991)
Algebraic-geometric codes.  Kluwer, Dordrecht
\endRefs
}
\bigskip

\noindent Russian Academy of Sciences

\noindent Steklov Mathematical Institute

\noindent St.~Petersburg Department

\noindent Fontanka 27, 

\noindent St.~Petersburg 191011, Russia

\noindent e-mail: skrig  \@ pdmi.ras.ru

\end